\documentclass[12pt]{amsart}

\usepackage{amsmath, amssymb, amsfonts, amsthm, amscd}
\usepackage{manfnt}
\usepackage{mathtools}
\usepackage{fullpage}
\usepackage{url}
\usepackage[all]{xy}
   \SelectTips{cm}{10}
\usepackage{txfonts}
\usepackage{enumitem}
\newtheorem{thm}{Theorem}[section]
\newtheorem{lemma}[thm]{Lemma}
\newtheorem{prop}[thm]{Proposition}

\newtheorem{hypothesis}[thm]{Hypothesis}
\newtheorem{prop-conj}[thm]{Proposition-Conjecture}

\theoremstyle{definition}
\newtheorem{defn}[thm]{Definition}

\theoremstyle{remark}
\newtheorem{rmk}[thm]{Remark}

\theoremstyle{remark}

\theoremstyle{remark}

\theoremstyle{remark}

\newcommand{\Q}{\mathbb{Q}}
\newcommand{\Qb}{\overline{\mathbb{Q}}}
\newcommand{\Z}{\mathbb{Z}}

\newcommand{\Ql}{\mathbb{Q}_{\ell}}
\newcommand{\Qlb}{\overline{\mathbb{Q}}_\ell}

\newcommand{\fl}{\mathbb{F}_{\ell}}
\newcommand{\bbP}{\mathbb{P}}

\DeclareMathOperator{\Hom}{Hom}

\DeclareMathOperator{\End}{End}
\DeclareMathOperator{\Aut}{Aut}
\DeclareMathOperator{\Out}{Out}

\DeclareMathOperator{\Sym}{Sym}
\DeclareMathOperator{\Ad}{Ad}
\DeclareMathOperator{\ad}{ad}

\DeclareMathOperator{\im}{im}

\DeclareMathOperator{\rk}{rk}

\DeclareMathOperator{\ord}{ord}

\newcommand{\gal}[1]{\Gamma_{#1}} 
\newcommand{\Gal}{\mathrm{Gal}} 

\newcommand{\onto}{\twoheadrightarrow}

\newcommand{\mc}{\mathcal}
\newcommand{\mf}{\mathfrak}
\newcommand{\mr}{\mathrm}
\newcommand{\mbf}{\mathbf}

\newcommand{\br}{\bar{\rho}} 
\newcommand{\fg}{\mathfrak{g}}
\newcommand{\fb}{\mathfrak{b}}

\DeclareMathOperator{\Lift}{\mathrm{Lift}}
\DeclareMathOperator{\Def}{\mathrm{Def}}

\newcommand{\fr}{\mathrm{fr}}
\newcommand{\Fr}{\mathrm{Fr}}
\newcommand{\kbar}{\overline{\kappa}}
\newcommand{\Ram}{\mathrm{Ram}}

\newcommand{\tv}{\tilde{v}}

\newcommand{\tF}{\widetilde{F}}
\newcommand{\LG}{{}^L G}

\begin{document}

\author{Stefan Patrikis}
\email{patrikis@math.utah.edu}
\address{The University of Utah \\ Salt Lake City, UT 84112}
\date{February 2015}

\title{Deformations of Galois representations and exceptional monodromy, II: raising the level}
\thanks{I am grateful to Chandrashekhar Khare and Bjorn Poonen for helpful conversations.}
\begin{abstract}
Building on lifting results of Ramakrishna, Khare and Ramakrishna proved a purely Galois-theoretic level-raising theorem for odd representations $\bar{\rho} \colon \Gal(\Qb/\Q) \to \mathrm{GL}_2(\overline{\mathbb{F}}_{\ell})$. In this paper, we generalize these techniques from type A1 to general (semi-)simple groups. We then strengthen our previous results on constructing geometric Galois representations with exceptional monodromy groups, achieving such constructions for almost all $\ell$, rather than a density-one set, and achieving greater flexibility in the Hodge numbers of the lifts; the latter improvement requires the new level-raising result.
\end{abstract}

\maketitle

\section{Introduction}
This paper enhances the deformation-theoretic techniques of \cite{stp:exceptional} and strengthens the applications in that paper to the construction of geometric Galois representations with exceptional algebraic monodromy groups. The foundation of the deformation-theoretic method of \cite{stp:exceptional} is an ingenious idea of Ravi Ramakrishna (\cite{ramakrishna:lifting}, \cite{ramakrishna02}), which shows that most odd representations $\br \colon \Gal(\Qb/\Q) \to \mr{GL}_2(\overline{\mathbb{F}}_{\ell})$ admit geometric (in the sense of Fontaine-Mazur) lifts to characteristic zero:
\[
\xymatrix{
& \mr{GL}_2(\overline{\Z}_{\ell}) \ar[d] \\
\gal{\Q} \ar[r]_-{\br} \ar@{-->}[ur]^{\rho} & \mr{GL}_2(\overline{\mathbb{F}}_{\ell}).
}
\]
``Odd" here means that the image $\br(c)$ of complex conjugation is conjugate to $\begin{pmatrix} 1 & 0 \\ 0 & -1 \end{pmatrix}$; a famous conjecture of Serre--now a theorem of Khare-Wintenberger (\cite{khare-wintenberger:serre1}, \cite{khare-wintenberger:serre2}), building on many other deep developments--asserts that any odd $\br$ admits a modular lift $\rho$, and in particular a lift of the sort produced by Ramakrishna's theorem. Ramakrishna's method also has implications for even representations, but it will not produce geometric lifts in this setting.
 
In \cite{stp:exceptional} we showed that Ramakrishna's ideas can be extended to suitable ``odd" representations valued in quite general reductive groups (we will specialize the setting somewhat for convenience). Let us recall this more precisely. For any field $F$, we set $\gal{F}= \Gal(\overline{F}/F)$ for some algebraic closure $\overline{F}$, and we now take $F$ to be a totally real number field. Let $G$ be a simply-connected, almost-simple group over $F$, and let $\LG$ denote a Langlands L-group for $G$ (fixing a pinned based root datum), which we regard as a split reductive group scheme over $\Z$. We denote its identity component, the Langlands dual group, by $G^\vee$, and we let $\fg^\vee$ denote the Lie algebra of $G^\vee$. Let $k$ be a finite extension of $\fl$, and let $\mc{O}$ be the ring of Witt vectors of $k$. A continuous homomorphism $\br \colon \gal{F} \to \LG(k)$ is \textit{odd} if for all $v \vert \infty$, with corresponding conjugacy classes of complex conjugation $c_v \in \gal{F}$, the local invariants satisfy
\[
\dim \left(\br(\fg^\vee)^{\Ad(\br(c_v))=1} \right)= \dim(fl_G),
\]
where $\br(\fg^\vee)$ is the adjoint representation, regarded as a $\gal{F}$-representation, and $fl_G$ is the flag variety of $G$. Not all such L-groups contain order-two elements $\br(c_v)$ satisfying this condition; here we note simply that if $-1$ belongs to the Weyl group $W_G$ of $G$, then $G^\vee$ always contains such ``odd'' elements, whereas if $-1$ is not in $W_G$, then such odd elements only exist in an L-group of the form $\LG= G^\vee \rtimes \Gal(\tF/F)$ where $\tF$ is a quadratic imaginary extension of $F$, and the outer automorphism class of $\Gal(\tF/F)$ acting on $G^\vee$ is a non-trivial $\Z/2$-symmetry of the Dynkin diagram of $G$. (For details, see \cite[\S 4.5, \S 9.1, \S 10.1]{stp:exceptional}.) The main deformation-theoretic results of \cite{stp:exceptional} (see \cite[Theorems 6.4, 7.4, 10.3, 10.4]{stp:exceptional}) show that odd homomorphisms $\br \colon \gal{F} \to \LG(k)$ satisfying appropriate global image and local ramification (especially at $\ell$) hypotheses always admit geometric lifts:
\[
\xymatrix{
& \LG(\mc{O}) \ar[d] \\
\gal{F} \ar[r]_-{\br} \ar@{-->}[ur]^{\rho} & \LG(k).
}
\]

Building on this deformation-theoretic machinery, we found the following application to the construction of geometric Galois representations with exceptional monodromy groups:
\begin{thm}[Theorems 8.4 and 10.6 of \cite{stp:exceptional}]\label{oldmono}
Let $G$ be a simply-connected group of exceptional type. Then there is a density-one set $\mc{L}$ of primes $\ell$, and for all $\ell \in \mc{L}$ a geometric Galois representation $\rho_{\ell} \colon \gal{\Q} \to \LG(\Qlb)$ with Zariski-dense image.  
\end{thm}
We note that for exceptional $G$ other than $\mr{E}_6$, we can always take $\LG= G^\vee$ in this theorem; for $G= \mr{E}_6$, we must take $\LG= G^\vee \rtimes \Out(\mr{E}_6)$. For $G$ of types $\mr{F}_4$ or $\mr{E}_6$, Theorem \ref{oldmono} was the first such construction of Galois representations over number fields with these monodromy groups. Beautiful work of Dettweiler-Reiter (for $G= \mr{G}_2$: \cite{dettweiler-reiter:rigidG2}) and especially Zhiwei Yun (for $G= \mr{G}_2, \mr{E}_7, \mr{E}_8$: \cite{yun:exceptional}) had previously established the other cases (indeed with stronger conclusions, that these $\rho_{\ell}$ could be found in the cohomology of algebraic varieties). We might ask for various strengthenings of Theorem \ref{oldmono}, but the first unsatisfactory aspect is that it only works for a density-one set of primes $\ell$; these arise as the set of ordinary primes of some well-chosen modular form.

As with \cite{stp:exceptional}, the present paper has two aims: the first is to extend to general groups a deformation-theoretic technique of Khare-Ramakrishna (\cite{khare-ramakrishna}) that strengthens Ramakrishna's original work, and the second is to apply this technique, and some other ideas, to strengthen the exceptional monodromy application in Theorem \ref{oldmono}. First we recall the results of \cite{khare-ramakrishna}, which continues in Ramakrishna's original setting of two-dimensional (let us say odd) $\br \colon \gal{\Q} \to \mr{GL}_2(k)$. Ramakrishna's method proceeds for $G= \mr{GL}_2$ by allowing into the level of lifts $\rho$ of $\br$ additional ramification of ``Steinberg-type'' at some carefully-chosen finite set $Q$ of primes. A curious circumstance results, however, in which the ultimate geometric lift $\rho$ is not known to be ramified at these auxiliary primes: roughly speaking, the Steinberg component of the local deformation ring will intersect the unramified component where the monodromy operator degenerates, and Ramakrishna's method cannot tell whether $\rho|_{\gal{\Q_q}}$ (for $q \in Q$) is at this intersection, or is a more general point of the Steinberg deformation ring. Khare and Ramakrishna found a more elaborate deformation-theoretic argument that allows one to ``force'' ramification at these auxiliary primes $q \in Q$ (precisely: to replace the initial set $Q$ with another auxiliary set where one can prove ramification). The reader would do well to keep in mind the automorphic and motivic analogues: establishing ramification of $\rho|_{\gal{\Q_q}}$ is the Galois-theoretic analogue of, respectively, the Ramanujan and weight-monodromy conjectures at the prime $q$. It should therefore not be surprising that there is something to prove!

Both the Ramakrishna method and the Khare-Ramakrishna method are highly sensitive to the image of $\br$. The construction of \cite{khare-ramakrishna} works for $\br(\gal{\Q}) \supset \mr{SL}_2(\mathbb{F}_{\ell})$, and in general the deformation-theoretic machinery functions most smoothly when $\br(\gal{\Q})$ contains the image of the $\fl$-points of the simply-connected cover of $G^\vee$. Another situation amenable to analysis, and essential in \cite{stp:exceptional}, arises when $\br(\gal{\Q})$ contains $\varphi(\mr{SL}_2(\fl))$, where 
\[
\varphi \colon \mr{SL}_2 \to G^\vee
\]
is a \textit{principal} homomorphism, i.e. the Jacobson-Morosov $\mr{SL}_2$ associated to a regular nilpotent element of $\fg^\vee$. In this paper we extend the ramification-forcing (or ``level-raising'') techniques of \cite{khare-ramakrishna} to both of these image settings: see Theorem \ref{raising} and Theorem \ref{sl2raising}, which are the technical heart of the paper. We use these along with some other deformation-theoretic work to deduce the following exceptional monodromy application (here $h_{G^\vee}$ denotes the Coxeter number of $G^\vee$):
\begin{thm}\label{mainintro}
Let $G$ be a simply-connected group of exceptional type. Then for all $\ell > 4h_{G^\vee}-1$, and in the case $G = \mr{E}_8$ also excluding $\ell= 229, 269, 367$, there are infinitely many geometric representations $\rho_{\ell} \colon \gal{\Q} \to \LG(\Ql)$ having Zariski-dense image. 

For all $G$ not of type $\mr{E}_6$, we can moreover arrange that the Hodge-Tate co-character (modulo conjugation) of $\rho_{\ell}$ is the half-sum of the positive co-roots $\rho^\vee$ of $G^\vee$.
\end{thm}
Note that in this theorem, the density one restriction on $\ell$ of Theorem \ref{oldmono} has been removed. The key to doing this is, instead of working with a fixed modular form and varying $\ell$, for each $\ell$ to construct a different elliptic curve whose mod $\ell$ representation $\bar{r}_{E, \ell}$ can be used as the seed residual representation; we then deform $\varphi \circ \bar{r}_{E, \ell}= \br \colon \gal{\Q} \to \LG(\fl)$. As the proof of Theorem \ref{excapp} will show, there is a great deal of flexibility in choosing these elliptic curves $E$, so we certainly get infinitely many examples in each case of Theorem \ref{mainintro}, all of which are different from the examples of Theorem \ref{oldmono}. We also get lifts with $\Ql$, rather than $\Qlb$, coefficients, in contrast to Theorem \ref{oldmono}. Moreover, we have achieved greater flexibility in the Hodge numbers of the lifts. We note especially the case of Hodge co-character $\rho^\vee$, since such lifts $\rho$ could conceivably, like the exceptional monodromy examples of \cite{yun:exceptional}, arise as specializations of an arithmetic local system over a curve. It would be interesting to see whether these examples could be used to reverse-engineer new motivic examples following Yun's techniques. Finally, we note that the original application of the level-raising method in \cite{khare-ramakrishna} was to produce two-dimensional $\rho \colon \gal{\Q} \to \mr{GL}_2(\mc{O})$ for which one could prove finiteness of a Selmer group associated to $\rho$ (see \cite[Theorem 4]{khare-ramakrishna}). Questions of this sort will be studied in the UCLA thesis of Mohammed Zuhair.
\section{Review of Ramakrishna's method}\label{ramreview}
We recall the set-up from \cite{stp:exceptional}. This will allow us to review the results we rely on, and to fix some important notation for the rest of the paper. If $\Psi$ is a based root datum for a simply-connected, almost-simple group, we fix a pinned split reductive group scheme $G^\vee$ over $\Z$ whose based root datum is dual to $\Psi$ (see \cite[\S 9.1]{stp:exceptional}. We then consider two cases:
\begin{itemize}
\item If $-1$ belongs to the Weyl group of $\Psi$, then we set ${}^L G= G^\vee$.
\item If $-1$ does not belong to the Weyl group of $\Psi$, then we let $\tF/F$ be a quadratic CM (totally imaginary) extension, and let $G$ be an $F$-form of the split group with root datum $\Psi$ such that the associated $\gal{F}$-action on $\Psi$ is non-trivial and factors through $\Gal(\tF/F)$. We then work with the associated L-group $\LG= G^\vee \rtimes \Gal(\tF/F)$; see \cite[\S 10.1]{stp:exceptional} for details.
\end{itemize}
The chosen pinning specifies a principal homomorphism $\varphi \colon \mr{PGL}_2 \to G^\vee$, defined over $\Z_{\ell}$ for all $\ell \geq h_{G^\vee}$ (see \cite[\S 2.4]{serre:principalsl2}), which extends to an L-homomorphism $\varphi \colon \mr{PGL}_2 \times \gal{F} \to \LG$. We recall (see \cite[Lemma 10.1]{stp:exceptional}) that 
\[
\theta_v= \varphi \left( \begin{pmatrix} 1 & 0 \\ 0 & -1 \end{pmatrix}, c_v \right)
\]
is a split Cartan involution for all complex conjugations $c_v \in \gal{F}$, i.e. $\dim (\fg^\vee)^{\Ad(\theta_v)}= \dim fl_G$. In particular, if $\bar{r}$ is an odd two-dimensional representation, then $\varphi \circ \bar{r}$ will be an odd $\LG$-valued representation.

For the time being (we will impose additional constraints as necessary), $\ell$ will be a prime that is at least 3 and is ``very good'' for $G$ (see \cite[\S 3]{stp:exceptional} for the implications), and let $k$ be a finite extension of $\fl$, with ring of Witt vectors $\mc{O}= W(k)$. Let $\Sigma$ be a finite set of places of $F$, assumed split in $\tF/F$ (if we are in the case $-1 \in W_G$, we set $\tF=F$), and we set $\gal{\Sigma}= \Gal(\tF(\Sigma)/F)$, where $\tF(\Sigma)$ is the maximal extension of $\tF$ in $\overline{F}$ that is unramified outside the places of $\tF$ above $\Sigma$ (for more notational details, see \cite[\S 9.2]{stp:exceptional}). In the generalization of Ramakrishna's method discussed in \cite{stp:exceptional}, we begin with an odd L-homomorphism
\[
 \br \colon \gal{\Sigma} \to \LG(k),
\]
subject to certain local (ramification) and global (image) hypotheses (see Hypothesis \ref{hyp0} and Theorem \ref{sl2raising} below), and we produce a geometric lift
\[
 \xymatrix{
& \LG(\mc{O}) \ar[d] \\
\gal{\Sigma \cup Q} \ar[r]_-{\br} \ar@{-->}[ur]^{\rho^Q} & \LG(k)
}
\]
for some auxiliary set of primes $Q$ of $F$, disjoint from $\Sigma$ and split in $\tF/F$. More precisely, under the just-mentioned hypotheses on $\br$, we define local deformation conditions $\mc{P}_v$ at each $v \in \Sigma$ (these are in fact defined by choosing one of the places $\tilde{v}$ of $\tF$ above $v$, and defining a local deformation condition for $\br|_{\gal{\tF_{\tilde{v}}}}$; the details of this don't concern us here, but see \cite[\S 9.2]{stp:exceptional}), which we abbreviate by $\mc{P}= \{\mc{P}_v\}_{v \in \Sigma}$, and Ramakrishna-type local deformation conditions $\mc{P}_q^{\mr{Ram}}$ (with respect to a fixed root of the maximal torus centralizing $\br(fr_q)$, but which we omit from the notation; see \cite[\S 4.2]{stp:exceptional}) at the additional primes in $Q$. We abbreviate the collection of all of these local conditions by $\mc{P}_Q= \{\mc{P}_v\}_{v \in \Sigma} \cup \{\mc{P}_q^{\mr{Ram}}\}_{q \in Q}$; this notation will be convenient because in our arguments we will have occasion to fix $\mc{P}$ and vary $Q$. Associated to these local conditions on $\br$ is a global deformation functor satisfying the corresponding local deformation conditions; it is representable, and our previous lifting theorems (see \cite[Theorems 6.4, 7.4, 10.3, 10.4]{stp:exceptional}; for more background on the deformation theory, we refer the reader to \cite[\S 3, \S 9.2]{stp:exceptional}) show that, for appropriately-chosen $Q$, the universal deformation ring $R_{\br}^{\mc{P}_Q}$ is isomorphic to $\mc{O}$. We will enshrine the terms of this conclusion in a definition:
\begin{defn}\label{auxset}
 We call a set of primes $Q$, distinct from $\Sigma$ and split in $\tF/F$, an \textit{auxiliary set} if for some choice of Ramakrishna-type deformation condition at primes in $Q$, the associated universal deformation ring $R_{\br}^{\mc{P}_Q}$ is isomorphic to $\mc{O}$.
\end{defn}
The lift $\rho^Q$ arises as some representative of the universal deformation; we denote its reduction modulo $\ell^n$ by $\rho^Q_n$.

Finally, as many of our calculations will involve manipulating Selmer groups with slightly different sets of local conditions, we fix some notation to avoid excessive clutter later on. For any set $T$ of primes of $F$ split in $\tF/F$ such that $\br$ defines a homomorphism $\gal{T} \to \LG(k)$, and for which we have specified a choice of extensions $\widetilde{T}= \{\tilde{w}\}_{w \in T}$ of the elements of $T$ to $\tF$, and for any set $\mc{L}=\{L_{\tilde{w}}\}_{w \in T}$ of subspaces $L_{\tilde{w}} \subset H^1(\gal{\tF_{\tilde{w}}}, \br(\fg^\vee))$, we set
\[
 H^1_{\mc{L}}= H^1_{\mc{L}}(\gal{T}, \br(\fg^\vee))= \ker \left(H^1(\gal{T}, \br(\fg^\vee)) \to \bigoplus_{w \in T} H^1(\gal{\tF_{\tilde{w}}}, \br(\fg^\vee))/L_{\tilde{w}} \right).
\]
 We will also abbreviate $h^1_{\mc{L}}= \dim_k(H^1_{\mc{L}})$. We let $L_{\tilde{w}}^{\perp}$ denote the orthogonal complement of $L_{\tilde{w}}$ under local duality, and we abbreviate the dual Selmer group associated to these local conditions by $H^1_{\mc{L}^\perp}$ (the Galois module is in this case $\br(\fg^\vee)(1)$). When the local subspaces $L_{\tilde{w}}$ are associated to a choice of deformation condition $\mc{P}_w$, then we will instead write $H^1_{\mc{P}}$, $h^1_{\mc{P}}$, etc.
\section{Forcing ramification}
The basic method described in \S \ref{ramreview} for adding auxiliary primes of ramification to kill a dual Selmer group does not \textit{a priori} produce lifts that are ramified at the auxiliary primes. An elaboration due to Khare and Ramakrishna (see \cite[Theorem 3]{khare-ramakrishna}) of Ramakrishna's method for the group $G= \mr{GL}_2$ does allow one to force ramification at the auxiliary primes. In the present section, we explain how to generalize the technique of \cite{khare-ramakrishna} to general $G$. 
\subsection{Axiomatics}
We begin with an `axiomatized' version of the argument for forcing ramification. Resuming the notation of \S \ref{ramreview}, we assume we have a residual representation $\br \colon \gal{\Sigma} \to {}^L G(k)$. To ensure that Ramakrishna's method will apply to $\br$, we assume it satisfies \cite[Properties (1)-(6) of \S 10.2]{stp:exceptional}; these properties encode an axiomatization of the lifting method, and all the reader needs to know is that they imply (\cite[Proposition 10.2]{stp:exceptional}) the existence of an auxiliary set of places (in the sense of Definition \ref{auxset}) $Q$ of $F$, disjoint from $\Sigma$ and split in $\tF/F$, and a lift $\rho^Q$ of type $\mc{P}_v$ for all $v \in \Sigma$ and Ramakrishna-type $\mc{P}_q^{\mr{Ram}}$ at all $q \in Q$,
\[
\xymatrix{
& {}^L G(\mc{O}) \ar[d] \\
\gal{\Sigma \cup Q} \ar[r]^{\br} \ar@{-->}[ur]^{\rho^Q} & {}^L G(k).
}
\]
Moreover, there are infinitely many choices of such auxiliary sets $Q$. In the following discussion, we will consider various characteristic zero lifts of $\br$, corresponding to different choices of auxiliary sets; recall that we systematically use the super-script notation $\rho^Q$ to indicate the unique deformation corresponding to the auxiliary set $Q$. Note that the deformation problem (as mentioned in \S \ref{ramreview}, and fully discussed in \cite[\S 9.2]{stp:exceptional}) has required fixing an extension of each prime $q \in Q$ to a place of $\tF$; to avoid complicating the notation, we will continue to denote by $q$ one such fixed extension.

By assumption, for the auxiliary set $Q$ we have the Selmer vanishing $h^1_{\mc{P}_Q}= h^1_{\mc{P}_Q^\perp}=0$. We will assume that $Q$ is non-empty, i.e. that the original deformation problem (of type $\mc{P}$) is obstructed; in \S \ref{unobs} we explain how to `raise the level' when $h^1_{\mc{P}^\perp}=0$.

We partition $Q$ as $Q_{\mr{ram}} \sqcup Q_{\mr{unr}}$, where a place $q \in Q$ belongs to $Q_{\mr{ram}}$ if and only if $\rho^Q$ is ramified at $q$. Our goal is to replace $Q$ with an auxiliary set for which $Q= Q_{\mr{ram}}$. We may assume that $Q_{\mr{unr}}$ is non-empty but chosen so that removing any prime $q$ from $Q_{\mr{unr}}$ yields a non-zero dual Selmer group ($h^1_{\mc{P}_{Q \setminus q}^\perp} \neq 0$). 
\begin{lemma}
Let $q$ be any element of $Q_{\mr{unr}}$, so that by assumption $Q_0= Q \setminus q$ is not an auxiliary set. Then $h^1_{\mc{P}_{Q_0}}= h^1_{\mc{P}_{Q_0}^\perp} =1$.
\end{lemma}
\proof
Following the notation of \cite[Proposition 5.2]{stp:exceptional}, let $L_q= L_q^{\mr{un}} \cap L_q^{\mr{Ram}}$;\footnote{We remind the reader precisely what this means with our notational conventions: $q$ is a place of $F$ with a fixed extension, also denoted $q$, to $\tF$, and $\br|_{\gal{\tF_q}}$ satisfies a Ramakrishna-type deformation condition with respect to some root $\alpha$, which is now fixed. The tangent space of this deformation condition is denoted $L_q^{\mr{Ram}}$, and the corresponding tangent space for unramified deformations is denoted $L_q^{\mr{un}}$.} this notation will be in effect for all of our arguments with auxiliary primes. The exact sequence
\[
0 \to H^1_{\mc{P}_Q^\perp} \to H^1_{\mc{P}_{Q_0}^\perp \cup L_q^\perp} \to L_q^{\perp}/L_q^{\mr{Ram}, \perp}
\]
implies ($Q$ is auxiliary) $h^1_{\mc{P}_{Q_0}^{\perp}} \leq h^1_{\mc{P}_{Q_0}^\perp \cup L_q^\perp} \leq 1$. The lemma follows since $h^1_{\mc{P}_{Q_0}^\perp}$ is non-zero by assumption.
\endproof
Fixing $q \in Q_{\mr{unr}}$ and letting $Q_0= Q \setminus q$, we can therefore choose bases $H^1_{\mc{P}_{Q_0}}= \langle \psi \rangle$ and $H^1_{\mc{P}_{Q_0}^\perp}= \langle \phi \rangle$. We now fix an integer $n$ such that $\rho^Q_{n-1}$ is ramified at every $v \in Q_{\mr{ram}}$. Forcing ramification depends on choosing two new auxiliary primes $q_1$ and $q_2$ satisfying the following criteria. We will assume for the time being that these can be arrange, and show how to find such primes in \S \ref{findaux} and \S \ref{sl2section}. 
\begin{hypothesis}\label{hyp1}
If $Q_{\mr{unr}}$ is non-empty, we assume there is a prime $q_1$, disjoint from $\Sigma \cup Q$ and split in $\tF/F$, satisfying (we continue to denote by $q_1$ a fixed extension to $\tF$)
\begin{itemize}
\item $\psi|_{\gal{\tF_{q_1}}} \in L_{q_1}$ (recall $L_{q_1}= L_{q_1}^{\mr{un}} \cap L_{q_1}^{\mr{Ram}}$);
\item $\phi|_{\gal{\tF_{q_1}}} \not \in L_{q_1}^{\mr{Ram}, \perp}$;
\item $\rho_{n-1}^Q|_{\gal{\tF_{q_1}}}$ is of Ramakrishna type, but $\rho_n^Q|_{\gal{\tF_{q_1}}}$ is not of Ramakrishna type (both taken with respect to the fixed root used in defining the Ramakrishna deformation condition at $q_1$)
\end{itemize}
\end{hypothesis}
Before proceeding to the conditions on $q_2$, we explain the first consequences of Hypothesis \ref{hyp1}:
\begin{lemma}\label{hypconseq}
Assume Hypothesis \ref{hyp1} holds (only the first two bulleted items are needed). Then
\begin{enumerate}
\item 
\[
H^1_{\mc{P}_{Q_0} \cup [L_{q_1}^{\mr{un}}+ L_{q_1}^{\mr{Ram}}]}= H^1_{\mc{P}_{Q_0 \cup q_1}} = H^1_{\mc{P}_{Q_0} \cup [L_{q_1}^{\mr{un}} \cap L_{q_1}^{\mr{Ram}}]}= \langle \psi \rangle.
\]
\item
$H^1_{\mc{P}_{Q_0 \cup q_1}^\perp}$ is one-dimensional, with a generator $\tilde{\phi}$, and $\phi$ and $\tilde{\phi}$ are linearly independent inside $H^1_{\mc{P}_{Q_0}^\perp \cup [L_{q_1}^\perp]}$.
\end{enumerate}
\end{lemma}
\proof
The subspace $H^1_{\mc{P}_{Q_0}^\perp \cup [L_{q_1}^{\mr{Ram}, \perp} \cap L_{q_1}^{\mr{un}, \perp}]}$ of $H^1_{\mc{P}_{Q_0}^\perp}= \langle \phi \rangle$ does not contain $\phi$, hence equals zero. From Wiles's formula (see \cite[Proposition 9.2]{stp:exceptional}), we deduce that
\[
h^1_{\mc{P}_{Q_0} \cup [L_{q_1}^{\mr{un}} + L_{q_1}^{\mr{Ram}}]}= \dim \left( L_{q_1}^{\mr{un}}+ L_{q_1}^{\mr{Ram}}\right)- \dim L_{q_1}^{\mr{un}}= 1,
\]
and again by the conditions on $q_1$ we deduce the first claim. In particular, $H^1_{\mc{P}_{Q_0 \cup q_1}^\perp}$ is one-dimensional, say with a generator $\tilde{\phi}$. Since $\phi$ is not contained in $H^1_{\mc{P}_{Q_0 \cup q_1}^\perp}$, we deduce the second claim.
\endproof
\begin{hypothesis}\label{hyp2}
If $Q_{\mr{unr}}$ is non-empty, we assume there is a prime $q_2$, disjoint from $\Sigma \cup Q \cup \{q_1\}$ and split in $\tF/F$, satisfying
\begin{itemize}
\item $\psi|_{\gal{\tF_{q_2}}} \not \in L_{q_2}$;
\item $\phi|_{\gal{\tF_{q_2}}} \not \in L_{q_2}^{\mr{un}, \perp} \cap L_{q_2}^{\mr{Ram}, \perp}$;
\item $\tilde{\phi}|_{\gal{\tF_{q_2}}} \not \in L_{q_2}^{\mr{un}, \perp} \cap L_{q_2}^{\mr{Ram}, \perp}$;
\item $\rho^Q_{n-1}|_{\gal{\tF_{q_2}}}$ is of Ramakrishna type.
\end{itemize}
\end{hypothesis}
\begin{lemma}\label{hypconseq2}
Under Hypotheses \ref{hyp1} and \ref{hyp2}, the sets $Q_0 \cup q_2$ and $Q_0 \cup \{q_1, q_2\}$ are both auxiliary. Moreover, as in Lemma \ref{hypconseq}, $H^1_{\mc{P}_{Q_0} \cup [L_{q_2}^{\mr{un}}+ L_{q_2}^{\mr{Ram}}]}= \langle \psi \rangle$.
\end{lemma}
\proof
We have already seen that $H^1_{\mc{P}_{Q_0}}= H^1_{\mc{P}_{Q_0 \cup \{q_1\}}}= \langle \psi \rangle$, $H^1_{\mc{P}_{Q_0}^\perp}= \langle \phi \rangle$, and $H^1_{\mc{P}_{Q_0 \cup \{q_1\}}^\perp}= \langle \tilde{\phi} \rangle$, so the first part of the lemma follows directly from the proof of \cite[Proposition 5.2, Proposition 10.2]{stp:exceptional}. The second part follows just as in Lemma \ref{hypconseq}.
\endproof
In particular, the universal deformation rings $R_{\br}^{\mc{P}_Q}$ and $R_{\br}^{\mc{P}_{Q_0 \cup \{q_1, q_2\}}}$ are both isomorphic to $\mc{O}$, with corresponding (lifts representing) universal deformations $\rho^Q$ (as before) and $\rho^{Q_0 \cup \{q_1, q_2\}}$. We have now assembled all the ingredients to prove the main result of this subsection:
\begin{prop}\label{formalforcing}
The mod $\ell^n$ reduction $\rho_n^{Q_0 \cup \{q_1, q_2\}}$ is ramified at all primes in $Q_{\mr{ram}} \cup \{q_1, q_2\}$. Consequently, there exist a (possibly empty) auxiliary set $\widetilde{Q}$ and a lift 
\[
\xymatrix{
& {}^L G(\mc{O}) \ar[d] \\
\gal{\Sigma \cup \widetilde{Q}} \ar[r]^{\br} \ar@{-->}[ur]^{\rho^{\widetilde{Q}}} & {}^L G(k)
}
\]
such that $\rho^{\widetilde{Q}}$ is ramified (of Ramakrishna type) at every auxiliary place $\tilde{q} \in \widetilde{Q}$, and of type $\mc{P}$ at the places in $\Sigma$.
\end{prop}
\proof
By construction, $\rho_{n-1}^Q$ is of Ramakrishna type at the primes $q_1$ and $q_2$, so we have an equality of deformations (and we may assume of lifts) $\rho^Q_{n-1}= \rho_{n-1}^{Q_0 \cup \{q_1, q_2\}}$ (here we use that both universal deformation rings are isomorphic to $\mc{O}$). The proposition will follow by comparing the mod $\ell^n$ reductions of these universal deformations. Note that, by hypothesis, $\rho_n^Q|_{\gal{\tF_{q}}}$ is unramified, and by construction of $q_1$, $\rho_n^Q|_{\gal{\tF_{q_1}}}$ is not of Ramakrishna type. The former observation implies we can write
\[
\rho_n^{Q_0 \cup \{q_1, q_2\}}= (1+ \ell^{n-1}h)\cdot \rho_n^Q
\]
for some non-zero cocycle $h \in H^1(\gal{\Sigma \cup Q_0 \cup \{q_1, q_2\}}, \br(\fg))$ ($h$ is unramified at $q$). Suppose that $\rho_n^{Q_0 \cup \{q_1, q_2\}}$ were unramified at $q_i$ for $i=1$ or $i=2$. Denote by $j$ the element of $\{1,2\} \setminus \{i\}$. In either case, we find that $h$ then belongs to $H^1_{\mc{P}_{Q_0} \cup [L_{q_{j}}^{\mr{un}}+L_{q_j}^{\mr{Ram}}]}$, which by Lemmas \ref{hypconseq} and \ref{hypconseq2} equals $\langle \psi \rangle$. But $\psi|_{\gal{\tF_{q_1}}}$ belongs to $L_{q_1}= L_{q_1}^{\mr{un}} \cap L_{q_1}^{\mr{Ram}}$, so $h|_{\gal{\tF_{q_1}}}$ belongs to $L_{q_1}^{\mr{Ram}}$, contradicting the fact that $\rho_n^{Q}|_{\gal{\tF_{q_1}}}$ is not of Ramakrishna type. 

Having established the first claim, the second claim of the Proposition follows inductively.
\endproof
\subsection{Finding auxiliary primes: maximal image}\label{findaux}
In this subsection we address, for $\br$ having maximal image, the heart of the problem of forcing ramification: finding auxiliary primes satisfying the Hypotheses \ref{hyp1} and \ref{hyp2}.

Recall that for simplicity we are assuming $G$ is simply-connected, so that $G^\vee$ is adjoint. We write $G_{\mr{sc}}^\vee$ for the simply-connected cover of $G^\vee$, and we will begin by treating the case of the simplest image hypothesis (compare \cite[\S 6]{stp:exceptional}), where there exists a subfield $k'$ of $k$ such that 
\[
\im \left( G_{\mr{sc}}^\vee(k') \to G^\vee(k') \right) \subseteq \br(\gal{\tF}) \subseteq G^\vee(k').
\]
For $G= \mr{SL}_2$ this hypothesis holds almost everywhere in the compatible system of mod $\ell$ representations associated to a non-CM classical modular form (\cite[Theorem 3.1]{ribet:modlimage2}). More precisely, the main result of this section (Theorem \ref{raising}) will be achieved under the hypotheses of \cite[Theorem 10.3]{stp:exceptional}; we recall them here, and assume they are in effect for the rest of the section:
\begin{hypothesis}\label{hyp0}
\begin{enumerate}
\item The degree $[\tF(\mu_{\ell}):\tF]$ is $\ell-1$.
\item There is a subfield $k' \subset k$ such that 
\[
\im \left( G^{\vee}_{\mr{sc}}(k') \to G^\vee(k') \right) \subset \br(\gal{\tF, \Sigma}) \subset G^\vee(k').
\]
\item $\ell -1$ is greater than the maximum of $8\cdot \# Z_{G^{\vee}_{\mr{sc}}}$ and
\[
\begin{cases}
\text{$(h_{G^\vee}-1)\# Z_{G^{\vee_{\mr{sc}}}}$ if $\# Z_{G^{\vee}_{\mr{sc}}}$ is even; or} \\
\text{$(2h_{G^\vee}-2) \# Z_{G^{\vee}_{\mr{sc}}}$ if $\# Z_{G^{\vee}_{\mr{sc}}}$ is odd.}
\end{cases}
\]
\item $\br$ is odd, i.e. for all complex conjugations $c_v$, $\Ad(\br(c_v))$ is a split Cartan involution of $\mf{g}^\vee$.
\item For all places $v \in \Sigma$ not dividing $\ell \cdot \infty$, $\br|_{\gal{\tF_{\tv}}}$ satisfies a liftable local deformation condition $\mc{P}_v$ with tangent space of dimension $h^0(\gal{F_v}, \br(\fg^\vee))$.\footnote{For $\ell$ sufficiently large relative to $\Sigma$, this condition will be shown always to hold in the forthcoming thesis of Jeremy Booher (\cite{booher:thesis}).}
\item For all places $v \vert \ell$, $\br|_{\gal{\tF_{\tv}}}$ is ordinary and satisfies the conditions (REG) and (REG*) of \cite[\S 4.1]{stp:exceptional}.
\end{enumerate}
\end{hypothesis}
The image hypothesis has the following implication (which is not the sharpest result, but is more than enough for our purposes):
\begin{lemma}\label{minfieldbigim}
 Assume $\ell > 3$, and the image of $\br$ satisfies
\[
 \im \left( G_{\mr{sc}}^\vee(\mathbb{F}) \to G^\vee(\mathbb{F}) \right) \subseteq \br(\gal{\tF}) \subseteq G^\vee(\mathbb{F})
\]
for some finite extension $\mathbb{F}$ of $\mathbb{F}_{\ell}$. Then $\mathbb{F}$ is the minimal field of definition of $\br(\fg^\vee)$ as $\gal{\Sigma}$ representation.
\end{lemma}
\proof
For $\mathbb{F}$ a finite field, and $r \colon \Gamma \to \Aut_{\mathbb{F}}(V)$ a semi-simple $\mathbb{F}$-representation of a finite group $\Gamma$, the minimum field of definition of $V$ is the extension of $\mathbb{F}_{\ell}$ generated by the coefficients of the characteristic polynomials of elements in the image of $r$ (\cite[Lemma 6.13]{deligne-serre}); this field is in turn (by Brauer-Nesbitt) the fixed field of the subgroup of automorphisms $\sigma$ of $\mathbb{F}$ such that $\sigma(V) \cong V$. Now let $V$ be the absolutely irreducible representation $\br(\fg^\vee)$. Let $b$ be a generator of the cyclic group $\mathbb{F}^\times$, and let $|\mathbb{F}|= \ell^f$ be the order of $\mathbb{F}$. We may assume $f \geq 2$, else there is nothing to show. By assumption, the image of $\br$ contains $\alpha^\vee(b)$ for any coroot $\alpha^\vee$ of $G^\vee$, and this element has eigenvalues on $\br(\fg^\vee)$ contained in the set $\{b^{\pm 3}, b^{\pm 2}, b^{\pm 1}, 1\}$; the containment may be strict, but the eigenvalue set at least contains $b^{\pm 2}$. If a non-trivial automorphism $\sigma$ of $\mathbb{F}$, which me may assume is $x \mapsto x^{\ell^{f-1}}$, fixes $V$, then it preserves the eigenvalue set of $\alpha(b)$, and in particular we must have $b^{2 \ell^{f-1}} \in \{b^{\pm 3}, b^{\pm 2}, b^{\pm 1}, 1\}$. As the order of $b$ is $\ell^f -1$, we deduce that $\ell^f-1 \leq 2 \ell^{f-1}+3$, hence $\ell \leq 2+ \frac{4}{\ell}$. The lemma follows.
\endproof
We will therefore replace $k$ by $k'$ (in the notation of Hypothesis \ref{hyp0}) in all that follows: nothing about the formal setup changes, but we will now be able to invoke the following lemma, which is elementary representation theory combined with the fact that the Brauer groups of finite fields are trivial.
\begin{lemma}[Lemma 7 of \cite{ramakrishna:lifting}]\label{minfield}
Let $k$ be a finite field, and let $\Gamma$ be any group. Suppose $\tau \colon \Gamma \to \mr{GL}_n(k)$ is an absolutely irreducible representation with minimal field of definition $k$. Then any non-zero $\mathbb{F}_{\ell}[\Gamma]$-submodule of $\tau$ is equal to all of $\tau$.
\end{lemma}
Next, if $\rho$ (eg, take $\rho$ equal to our $\rho^Q$) is a lift to $\LG(\mc{O})$ of $\br$, we will need to understand what the image hypothesis implies about the image $\rho_n(\gal{\tF})$, so we continue with some group theory. 
\begin{lemma}\label{grouptheory}
Assume $\ell > 5$. Let $H$ be a split reductive group over $\mc{O}=W(k)$, and let $\mf{h}$ denote the Lie algebra of $H$ (over $\mc{O}$). Suppose that $H_1$ is a subgroup of $H(k)$ such that
\begin{enumerate}
\item $H_1$ acts absolutely irreducibly on $\mf{h}_k= \mf{h} \otimes_{\mc{O}} k$, with minimal field of definition equal to $k$.
\item Letting $\theta$ denote the highest root of $\mf{h}$ and $X_{\theta}$ be an $mc{O}$-basis of the root space $\mf{h}_{\theta}$, $H_1$ contains $\exp(X_{\theta})$ (here we write $\exp \colon \mf{h}_{\theta} \to H$ for the homomorphism giving the root subgroup corresponding to $\theta$, as in \cite[Theorem 4.1.4]{conrad:luminy}).
\end{enumerate}
Then any subgroup $H_n$ of $H(\mc{O}/\ell^n)$ whose image mod $\ell$ contains $H_1$ must contain 
\[
\ker \left(H(\mc{O}/\ell^n) \to H(\mc{O}/\ell) \right).
\]
In particular, if $H_1$ contains the image of $H^{\mr{sc}}(k) \to H(k)$, then $H_n= H(\mc{O}/\ell^n)$.
\end{lemma}
\begin{rmk}
The first assumption implies the center of $H$ is trivial, which will be the case in our application; but the statement and argument admit straightforward modifications in the case where $H_1$ acts absolutely irreducibly on the derived subalgebra of $\mf{h}_k$.
\end{rmk}
\proof
We argue by induction. Assume the claim has been established mod $\ell^n$, and fix a subgroup $H_{n+1}$ of $H(\mc{O}/\ell^{n+1})$ as in the lemma; we must show that the surjection $H_{n+1} \to H_n$ ($H_n$ being by definition the image mod $\ell^n$) has kernel equal to $\ker(H(\mc{O}/\ell^{n+1}) \to H(\mc{O}/\ell^n))$. The essential step is ruling out the existence of a section $H_n \to H_{n+1}$ of the reduction map. Suppose such a section existed, and let $u \in H_{n+1}$ be the image under this section of $\exp(X_{\theta})$. Then $u^{\ell^n}=1$, and $\Ad(u)^{\ell^n}=1$. Since $\ell \neq 2$, we find (doing the calculation in $\mf{h}$ and reducing)
\[
\Ad(u)= (1+ \ad(X_{\theta})+ \frac{\ad(X_{\theta})^2}{2}+ \ell^n X)
\]
for some $X \in \End_{\mc{O}}(\mf{h})$ (we have used that $\theta$ is the highest root to conclude $\ad(X_\theta)^3=0$). Now we expand
\begin{align*}
1= \Ad(u)^{\ell^n}&= \sum_{j=0}^{\ell^n} \binom{\ell^n}{j} \left(\ad(X_{\theta})+ \frac{\ad(X_{\theta})^2}{2}+ \ell^n X \right)^j \\
&\equiv 1+ \ell^n\left( \ad(X_{\theta})+ \frac{\ad(X_{\theta})^2}{2}\right)+ \binom{\ell^n}{2}\ad(X_{\theta})^2+ \left(\ad(X_{\theta})+\frac{\ad(X_{\theta})^2}{2}+ \ell^nX \right)^{\ell^n} \pmod {\ell^{n+1}}.
\end{align*}
The only surviving terms from this last binomial expansion have the form
\[
\left(\ad(X_{\theta})+ \frac{\ad(X_{\theta})^2}{2} \right)^i \ell^n X \left(\ad(X_{\theta})+ \frac{\ad(X_{\theta})^2}{2} \right)^{\ell^n-1-i}
\]
(using $\ad(X_{\theta})^{\ell^n-1}=0$ for $\ell \geq 5$), and these terms vanish unless $i \leq 2$ and $\ell^n-1-i\leq 2$; in particular, they vanish unless $\ell \leq 5$. We conclude that for $\ell >5$ (our running hypothesis), 
\[
\Ad(u)^{\ell^n} \equiv 1+ \ell^n \ad(X_{\theta}) \pmod{\ell^{n+1}},
\]
and thus $u$ cannot have order $\ell^n$, contradicting the assumption that there is a section $H_n \to H_{n+1}$.

In particular, $\ker \left(H_{n+1} \to H_n \right)$ is a non-trivial subgroup of $H_{n+1}$, necessarily stable under the adjoint action of $H_n$ (factoring through $H_1 \subset H(k)$) on 
\[
\mf{h}_k \otimes_k \ell^n \mc{O}/\ell^{n+1} \mc{O} \xrightarrow[\exp]{\sim} \ker \left(H(\mc{O}/\ell^{n+1}) \to H(\mc{O}/\ell^n) \right)
\]
Recall we have also assumed that $\mf{h}_k$ is an absolutely irreducible $k[H_1]$-module. Since the Brauer group of a finite field is trivial, in fact any non-zero $\mathbb{F}_{\ell}[H_1]$-submodule of $\mf{h}_k$ must then equal all of $\mf{h}_k$. Thus, the kernel of the reduction $H_{n+1} \to H_n$ must equal the full $\ker(H(\mc{O}/\ell^{n+1}) \to H(\mc{O}/\ell^n))$.

The application when $H_1$ contains $\im \left(H^{\mr{sc}}(k) \to H(k) \right)$ follows from Lemma \ref{minfieldbigim}.
\endproof
\begin{rmk}
Except for the case $\ell=5$, this generalizes (from the case $H= \mr{SL}_2$) \cite[IV.3.4 Lemma 3]{serre:ladic}. 
\end{rmk}
We now have the relevant group theory in place to construct auxiliary primes. Recall that we start with non-trivial Selmer classes $\psi$ and $\phi$ spanning $H^1_{\mc{P}_{Q_0}}$ and $H^1_{\mc{P}_{Q_0}^\perp}$, respectively. After restriction to $\gal{K}$, for $K= \tF(\br(\fg^\vee), \mu_{\ell})$, these become homomorphisms cutting out extensions $K_\psi$ and $K_{\phi}$, which are moreover Galois over $F$ (and $\tF$). Also let $P_n$ be the fixed field of $\rho^Q_n |_{\gal{K}}$ (note that $P_1=K$); it too is Galois over $F$. The essential Galois-theoretic point is the following:
\begin{lemma}\label{lindisjoint}
Assume $\ell > 5$ and that $\br(\gal{\tF, \Sigma})$ contains $\im \left(G^{\vee, \mr{sc}}(k) \to G^\vee(k) \right)$. Then (for all $n \geq 1$) the extensions $K_{\phi}$, $K_{\psi}$, $P_n$, and $K(\mu_{\ell^n})$ are strongly linearly disjoint over $K$, i.e. the intersection of any one with the compositum of the other three is equal to $K$.
\end{lemma}
\proof
We first show that $K(\mu_{\ell^n}) \cap P_n K_{\phi} K_{\psi}=K$. For future reference (see Proposition \ref{forceramsl2}), we observe that this part of the argument applies verbatim in the context of \S \ref{sl2section}.
The abelianization of the image of $\br$ has order prime to $\ell$ (see \cite[Lemma 6.6]{stp:exceptional}), so $K$ and $\tF(\mu_{\ell^n})$ are linearly disjoint over $\tF(\mu_{\ell})$, and therefore the conjugation action of $\Gal(K/\tF(\mu_{\ell}))$ on $\Gal(K(\mu_{\ell^n})/K)$ is trivial. Now, either there is some $i \leq n$ such that $P_i K_{\phi} K_{\psi} \cap K(\mu_{\ell^n})$ properly contains $P_{i-1} K_{\phi} K_{\psi} \cap K(\mu_{\ell^n})$; or $P_n K_{\phi} K_{\psi} \cap K(\mu_{\ell^n})= K_{\phi}K_{\psi} \cap K(\mu_{\ell^n})$. In the first case, we conclude that $\Gal(P_i K_{\phi} K_{\psi} \cap K(\mu_{\ell^n})/P_{i-1} K_{\phi} K_{\psi} \cap K(\mu_{\ell^n}))$ is, as $\mathbb{F}_{\ell}[\Gal(K/\tF)]$-module, a sub-quotient of $\br(\fg^\vee)$ isomorphic to a sum of copies of the trivial representation; and in the second, that $\Gal(K_{\phi}K_{\psi} \cap K(\mu_{\ell^n})/K)$ is a $\mathbb{F}_{\ell}[\Gal(K/\tF)]$-subquotient of $\br(\fg^\vee) \oplus \br(\fg^\vee)(1)$ isomorphic to a sum of copies of the trivial representation. (Note that the $\Gal(K/\tF)$ action on these modules factors through $\Gal(\tF(\mu_{\ell})/\tF)$ by the observation at the start of the proof about the abelianized image of $\br$; the remaining $\Gal(\tF(\mu_{\ell})/\tF)$-action is trivial because $\tF(\mu_{\ell^n})/\tF$ is abelian.) Our assumptions on $\br$ forbid this unless these subquotients are zero, and we conclude that $P_n K_{\phi} K_{\psi} \cap K(\mu_{\ell^n})= K$.

Next we show $P_n \cap K_{\phi} K_{\psi}=K$ for all $n \geq 1$. To do this, we inductively prove that $P_n \cap P_{n-1} K_{\phi} K_{\psi}= P_{n-1}$. Lemmas \ref{minfield} and \ref{grouptheory} imply that $\Gal(P_n/P_{n-1})$ is an irreducible $\mathbb{F}_{\ell}[\Gal(K/\tF)]$-module, so $\Gal(P_n \cap P_{n-1} K_{\phi} K_{\psi}/P_{n-1})$, being a $\Gal(P_{n-1}/\tF) \onto \Gal(K/\tF)$-stable quotient of $\Gal(P_n/P_{n-1})$, must be either trivial or all of $\Gal(P_n/P_{n-1})$. We cannot have $P_n \subset P_{n-1}K_{\phi}K_{\psi}$, however, since the extension
\[
1 \to \Gal(P_{n-1}K_{\phi}K_{\psi}/P_{n-1}) \to \Gal(P_{n-1} K_{\phi} K_{\psi}/\tF) \to \Gal(P_{n-1}/\tF) \to 1
\]
splits (the corresponding $H^2$ class is a sum of the images of $\phi$ and $\psi$ in the inflation-restriction-transgression sequence), whereas (by Lemma \ref{grouptheory}) the extension
\[
1 \to \Gal(P_n/P_{n-1}) \to \Gal(P_n/\tF) \to \Gal(P_{n-1}/\tF) \to 1
\]
does not split.

We have already shown (\cite[Lemma 6.8]{stp:exceptional}) that $K_{\phi} \cap K_{\psi}=K$, so we deduce as desired that
\[
\Gal(P_nK_\phi K_\psi (\mu_{\ell^n})/K) \xrightarrow{\sim} \Gal(P_n/K) \times \Gal(K_{\phi}/K) \times \Gal(K_{\psi}/K) \times \Gal(K(\mu_{\ell^n})/K).
\]
\endproof
Using this linear disjointness, we can arrange the criteria of Hypothesis \ref{hyp1} on auxiliary primes by working independently in these four Galois extensions of $K$:
\begin{prop}\label{findingprimes}
There exists a positive density of primes $q_1$ of $\tF$ such that $\phi|_{q_1}$, $\psi|_{q_1}$ and $\rho^Q_n|_{q_1}$ satisfy the criteria of Hypothesis \ref{hyp1}.
\end{prop}
\proof
Begin with $\sigma_1 \in \Gal(K/\tF)$ such that $\br(\sigma_1)$ is regular semi-simple and for some simple root $\alpha$ of the corresponding maximal torus (centralizing $\br(\sigma_1)$), $\alpha(\br(\sigma_1))= \kbar(\sigma_1)$: we can as in \cite[Lemma 6.7]{stp:exceptional} take $\br(\sigma_1)$ of the form $2\rho^\vee(t_1)$ for some $t_1 \in \mathbb{F}_{\ell}^\times$ of sufficiently large order ($2 \rho^\vee$ denotes the sum of the positive coroots of $G^\vee$). As usual, we fix from now on an $\alpha$ in defining deformations of Ramakrishna type. Choose an element $t_{n-1} \in (\mc{O}/\ell^{n-1})^\times$ lifting $t_1$ such that $t_{n-1}^2$ is in the image of $\kappa \colon \Gal(K(\mu_{\ell^{n-1}})/\tF) \to (\mc{O}/\ell^{n-1})^{\times}$; then Lemma \ref{grouptheory} and the disjointness statement $K(\mu_{\ell^{n-1}}) \cap P_{n-1}=K$ (Lemma \ref{lindisjoint}) imply we can find a lift $\sigma_{n-1}$ of $\sigma_1$ to $\Gal(K(\mu_{\ell^{n-1}})P_{n-1}/\tF)$ such that $\kappa(\sigma_{n-1})=t_{n-1}^2$ and $\rho^Q_{n-1}(\sigma_{n-1})= 2 \rho^\vee(t_{n-1})$. By the same reasoning, we can find lifts $t_n$ and $\sigma_n$ of $t_{n-1}$ and $\sigma_{n-1}$ such that $\kappa(\sigma_n) =t_n^2 \pmod {\ell^n}$ but $\rho_n^Q(\sigma_n)= 2 \rho^\vee(t_n+\ell^{n-1})$. Note that $\alpha(\rho_n^Q(\sigma_n)) \neq \kappa(\sigma_n) \pmod {\ell^n}$.

Now let $\sigma$ denote any extension of $\sigma_n$ to an element of $\Gal(P_n(\mu_{\ell^n})K_{\phi}K_{\psi}/\tF)$. By Lemma \ref{lindisjoint} and \cite[Lemma 6.7]{stp:exceptional} we can modify $\sigma$ by elements $\tau_{\phi}$ and $\tau_{\psi}$ of $\Gal(K_{\phi}/K)$ and $\Gal(K_{\psi}/K)$ (canonically lifted to $\Gal(P_{n}(\mu_{\ell^n})K_{\phi}K_{\psi}/K)$) so that $\phi(\tau_\phi \tau_\psi \sigma)$ has non-zero $\fg^\vee_{-\alpha}$ component, and so that $\psi(\tau_\phi \tau_\psi \sigma)$ is zero. The \v{C}ebotarev density theorem yields a positive density set of places $w$ of $\tF$ (we may assume split over $F$) whose frobenii $\fr_w$ lie in the conjugacy class of $\tau_\phi \tau_\psi \sigma$, and therefore satisfy (compare \cite[Lemma 5.3]{stp:exceptional})
\begin{itemize}
\item $\rho^Q_{n-1}(\fr_w)$ is of Ramakrishna type but $\rho^Q_n(\fr_w)$ is not of Ramakrishna type, with respect to our fixed root $\alpha$ (to be precise, to check the latter statement we must invoke the claim proved in the second paragraph of \cite[Lemma 4.10]{stp:exceptional});
\item $\phi|_{w} \not \in L_w^{\mr{Ram}, \perp}$;
\item $\psi|_w=0$ (and in particular belongs to $L_w= L_w^{\mr{un}} \cap L_w^{\Ram}$).
\end{itemize}
\endproof
In the same fashion, we can also achieve the conditions of Hypothesis \ref{hyp2}. Recall that once Hypothesis \ref{hyp1} was satisfied, we produced (see Lemma \ref{hypconseq}) a second dual Selmer element $\tilde{\phi} \in H^1_{\mc{P}_{Q_0}^\perp \cup [L_{q_1}^\perp]}$, linearly independent from $\phi$.
\begin{prop}\label{twolindisj}
There exists a positive density set of primes $q_2$ of $\tF$ such that $\psi|_{q_2}$, $\phi|_{q_2}$, $\tilde{\phi}|_{q_2}$, and $\rho_{n-1}^Q|_{q_2}$ satisfy the criteria of Hypothesis \ref{hyp2}.
\end{prop}
\proof
The argument reduces to arguments already given (Lemma \ref{lindisjoint}, Proposition \ref{findingprimes} and \cite[Lemma 5.3]{stp:exceptional}), provided we also establish linear disjointness of $K_{\phi}$ and $K_{\tilde{\phi}}$ over $K$. This too is rather standard: any intersection would give an $\mathbb{F}_{\ell}[\Gal(K/\tF)]$-stable quotient of $\Gal(K_{\phi}/K) \xrightarrow[\sim]{\phi} \br(\fg^\vee)(1)$, so that either $K_{\phi}=K_{\tilde{\phi}}$, or $K_{\phi} \cap K_{\tilde{\phi}}=K$. In the former case, the composite
\[
\br(\fg^\vee)(1) \xrightarrow{\phi^{-1}} \Gal(K_{\phi}/K) \xrightarrow{\tilde{\phi}} \br(\fg^\vee)(1)
\]
lies in $\End_{\mathbb{F}_{\ell}[\gal{\tF}]}(\br(\fg^\vee))$, which is $k$, since $k$ is the minimal field of definition of $\br(\fg^\vee)$. Thus, in this case $\phi$ and $\tilde{\phi}$ would be $k$-linearly dependent.
\endproof
Invoking Proposition \ref{formalforcing}, we obtain the main result of this section:
\begin{thm}\label{raising}
Let $\br \colon \gal{\Sigma} \to \LG(k)$ be a continuous representation satisfying Hypothesis \ref{hyp0}, i.e. the hypotheses of \cite[Theorem 10.3]{stp:exceptional}. Then there exist a finite set of primes $R$ of $F$, disjoint from $\Sigma$ and split in $\tF/F$, and a geometric lift $\rho \colon \gal{\Sigma \cup R} \to \LG(W(k))$ of $\br$, such that for each place $\tilde{v}$ of $\tF$ above an element of $R$, $\rho|_{\gal{\tF_{\tilde{v}}}}$ is ramified of Ramakrishna type.  
\end{thm}
\subsection{The case of an unobstructed initial deformation problem}\label{unobs}
Theorem \ref{raising} can be regarded as a (modest) level-raising result for the original lift $\rho^Q$ of $\br$ produced by \cite[Theorem 10.3]{stp:exceptional}. It may of course happen in that situation that we can take $Q= \emptyset$, in which case Theorem \ref{raising} is vacuous. We now explain how it is still possible to raise the level of such the lift $\rho^{\emptyset}$ that exists when $H^1_{\mc{P}}(\gal{\Sigma}, \br(\fg^\vee))=0$. I am grateful to Shekhar Khare for pointing this out.
\begin{prop}\label{unobstructed}
Suppose that $H^1_{\mc{P}}(\gal{\Sigma}, \br(\fg^\vee))=0$. Then there exists a prime $q$ of $\tF$ (split over a prime $q$ of $F$) such that $\br|_{\gal{\tF_{q}}}$ is of Ramakrishna type and $H^1_{\mc{P} \cup L_q^{\mr{Ram}}}$ is non-zero. Replacing $\Sigma$ by $\Sigma \cup \{q\}$, we can therefore find a non-empty finite set of primes $R$ as in Theorem \ref{raising} and a geometric lift $\rho \colon \gal{\Sigma \cup \{q\} \cup R} \to \LG(\mc{O})$ such that for each place $v$ of $R$ (with fixed extension $v$ to $\tF$), $\rho|_{\gal{\tF_{v}}}$ is ramified of Ramakrishna type.
\end{prop}
\proof
Denote by $\rho^{\emptyset}$ the $\ell$-adic lift produced by the hypothesis $H^1_{\mc{P}}=0$. We may assume that for all auxiliary primes $q$ such that $\br|_{\gal{\tF_{q}}}$ is of Ramakrishna type, $H^1_{\mc{P}\cup L_q^{\mr{Ram}}}=0$. As in the proof of Proposition \ref{findingprimes}, choose such a $q$ for which $\rho^\emptyset_2|_{\gal{\tF_{q}}}$ is \textit{not} of Ramakrishna type. By assumption, we obtain a unique $\ell$-adic deformation $\rho^{\{q\}}$ of type $\mc{P} \cup \mc{P}_q^{\mr{Ram}}$. If it is ramified at $q$, then we are done (take $R= \{q\}$). If not, then since each of the universal deformation rings $R^{\mc{P}}_{\br}$ and $R_{\br}^{\mc{P} \cup \mc{P}_q^{\mr{Ram}}}$ is simply $\mc{O}$, we must have an equality of deformations $\rho^{\{q\}}= \rho^{\emptyset}$. This contradicts the assumption that $\rho^{\emptyset}_2|_{\gal{\tF_{q}}}$ is not of Ramakrishna type, so we are done.
\endproof
\subsection{Finding auxiliary primes: the principal $\mr{SL}_2$}\label{sl2section}
With a view toward exceptional monodromy applications, we now turn to the case of $\br$ of the form $\varphi \circ \bar{r}$, where recall (see \S \ref{ramreview}) $\varphi \colon \mr{PGL}_2 \times \gal{F} \to \LG$ denotes the principal homomorphism, and $\bar{r}$ is a two-dimensional representation $\bar{r} \colon \gal{\Sigma} \to \mr{GL}_2(k) \times \gal{F}$. In particular, we now need $\ell \geq h_{G^\vee}$ in order to define the principal homomorphism $\ell$-integrally. Throughout, $F$ will be a totally real field satisfying $[F(\mu_{\ell}):F]= \ell -1$. Recall that in the process of defining $G^\vee$ we have fixed a pinning: letting $\Delta$ denote the simple roots of $G^\vee$, for all $\alpha \in \Delta$ we let $X_{\alpha}$ denote the pinned $\mc{O}$-basis of the root space $\fg^\vee_{\alpha}$. The principal homomorphism $\varphi$ satisfies
\[
d\varphi \begin{pmatrix} 0 & 1 \\ 0 & 0 \end{pmatrix}= \sum_{\alpha \in \Delta} X_{\alpha},
\]
and we set $X= \sum_{\alpha \in \Delta} X_{\alpha}$. We now prove the group theory lemma in this setting that will play the role of Lemma \ref{grouptheory}; note that this version is somewhat weaker.
\begin{lemma}\label{sl2grouptheory}
Assume $\ell > 4 h_{G^\vee}-3$, and that $\bar{r}(\gal{\tF})$ contains $\begin{pmatrix} 1 & 1 \\ 0 & 1 \end{pmatrix}$. Then there is no section $\bar{r}(\gal{\tF}) \to G^\vee(\mc{O}/ \ell^2)$.
\end{lemma}
\proof
If there were a section, mapping $\exp(X)= \varphi(\begin{pmatrix} 1 & 1 \\ 0 & 1 \end{pmatrix})$ to $u \in G^\vee(\mc{O}/ \ell^2)$, then
\[
1= \Ad(u)^\ell= [\exp(ad X) \exp( \ell Y)]^{\ell}
\]
for some $Y \in \End_k(\fg^\vee)$. The argument is similar to that of Lemma \ref{grouptheory}, but now using the fact that $ad(X)^{2 h_{G^\vee} -1}=0$. Namely, set $Z= \sum_{i \geq 1} \frac{\ad(X)^i}{i!}$, so that we are to compute 
\[
\left[(1+Z)(1+\ell Y)\right]^\ell= \sum_{j=0}^\ell \binom{\ell}{j}(Z+ \ell(Y+ZY))^j,
\]
which is equal modulo $\ell^2$ to $(1+Z)^\ell$ plus a sum of terms of the form $Z^j \ell(Y+ZY) Z^{\ell-1-j}$, the latter expression vanishing unless $j \leq 2h_{G^\vee}-2$ and $\ell-1-j \leq 2h_{G^\vee}-2$, and in particular unless $\ell \leq 4h_{G^\vee} -3$. Thus, under our hypotheses, we would have
\[
1= \Ad(u)^\ell= (1+Z)^\ell= \exp(\ad(\ell X)) \pmod{\ell^2},
\]
contradicting the fact that the $X_{\alpha}$ are bases of the free $\mc{O}$-modules $\fg^\vee_{\alpha}$.
\endproof
We now show how to find auxiliary primes in the setting of this section. That is, we now assume that the image of $\bar{r} \colon \gal{\Sigma} \to \mr{GL}_2(k)$ satisfies
\[
 \mr{SL}_2(k') \subset \bar{r}(\gal{\tF}) \subset k^\times \cdot \mr{GL}_2(k')
\]
 for some subfield $k'$ of $k$; since we will only work with the projectivization of $\bar{r}$, we replace $k$ by $k'$ in all that follows. Provided $\ell \geq 2h_{G^\vee}-1$, we find that (this new) $k$ is the minimal field of definition of each of the (absolutely) irreducible factors in the decomposition of the Lie algebra (see \cite[Lemma 7.3]{stp:exceptional}),
\[
\Ad(\varphi \circ \bar{r}) \cong \bigoplus_{i=1}^l \Sym^{2 m_i}(k^2)\otimes {\det}^{-m_i} (\bar{r}),
\]
by the following lemma (note that the maximum $m_i$ appearing is $h_{G^\vee} -1$):
\begin{lemma}
Let $r$ be a positive integer, and let $\mathbb{F}$ be a finite extension of $\mathbb{F}_{\ell}$. If $\ell > r+1$, then the minimal field of definition of $\mr{SL}_2(\mathbb{F})$ on $\Sym^r(\mathbb{F}^2)$ is $\mathbb{F}$.
\end{lemma}
\proof
As in Lemma \ref{minfieldbigim}, it suffices for a generator $b$ of $\mathbb{F}^\times$ to compare the eigenvalue sets $\{b^r, b^{r-2}, \ldots, b^{-r}\}$ and $\{b^{r \ell^{f-1}}, b^{(r-2)\ell^{f-1}}, \ldots, b^{-r \ell^{f-1}}\}$. If they agree, then the order of $b$ is bounded above by $r+r \ell^{f-1}$, i.e. $r \geq \frac{\ell^f-1}{\ell^{f-1}+1}$. For $r < \ell -1$ we obtain a contradiction.
\endproof
\begin{prop}\label{forceramsl2}
There is a positive density set of primes $q_1$ of $\tF$ satisfying the criteria of Hypothesis \ref{hyp1}. Having found such a prime $q_1$, and used it to produce a cohomology class $\tilde{\phi}$ as in Lemma \ref{hypconseq}, there is a positive density set of primes $q_2$ satisfying the criteria of Hypothesis \ref{hyp2}.
\end{prop}
\proof
Linear disjointness of $K(\mu_{\ell^2})$ from $P_2 K_\phi K_\psi$ follows from the same argument as in Lemma \ref{lindisjoint}, since $\br(\fg^\vee)$ and $\br(\fg^\vee)(1)$ contain no copy of the trivial representation. Now consider the intersection $P_2 \cap K_\phi K_\psi$; as in Lemma \ref{lindisjoint}, the extension
\[
1 \to \Gal(P_2 \cap K_\phi K_\psi /K) \to \Gal(P_2 \cap K_\phi K_\psi/\tF) \to \Gal(K/\tF) \to 1
\]
splits. We fix a section $s$, set $u= s(\varphi(\begin{pmatrix} 1& 1 \\ 0 & 1 \end{pmatrix}))$, and then let $\tilde{u}$ be some choice of lift of $u$ to an element of $\Gal(P_2/\tF) \subset G^\vee(\mc{O}/\ell^2)$. By the calculation in Lemma \ref{sl2grouptheory}, $\Ad(\tilde{u}^\ell)= \exp (\ad (\ell X)) \pmod {\ell^2}$ (recall $X$ is the regular nilpotent element produced from the pinning), so $\tilde{u}^\ell$ is non-trivial. At the same time, by construction $\tilde{u}^\ell$ lies in $\Gal(P_2/P_2 \cap K_\phi K_\psi)$. 
We deduce that $\Gal(P_2/P_2 \cap K_{\phi} K_{\psi})$ at least contains the constituent $(\Sym^2 \otimes {\det}^{-1})(\bar{r}) \otimes_k \ell \mc{O}/\ell^2 \mc{O}$ of $\ker(G^\vee(\mc{O}/\ell^2) \to G^\vee(\mc{O}/\ell))$. This will suffice for our purposes.

As in the proof of Proposition \ref{findingprimes}, we choose $\sigma_1 \in \Gal(K/\tF)$ such that $\br(\sigma_1)= 2 \rho^\vee(t_1)$ ($t_1 \in \mathbb{F}_{\ell}^\times$) has the desired properties modulo $\ell$. We fix an extension to $\Gal(P_2K_{\phi}K_{\psi}(\mu_{\ell^2})/\tF)$ and then use linear disjointness of $K_{\phi}$ and $K_{\psi}$ over $K$ to modify this to an element $\tilde{\sigma}_1$ such that $\psi(\tilde{\sigma}_1)=0$ and $\phi(\tilde{\sigma}_1)$ has non-zero $\fg^\vee_{-\alpha}$ component for the root $\alpha$ specified in our `Ramakrishna-type' deformation condition. Now, $\rho^Q_2(\tilde{\sigma}_1)$ may or may not be $\widehat{G^\vee}(\mc{O}/\ell^2)$-conjugate to $H_{\alpha}=U_{\alpha} \cdot T^\vee$ ($U_{\alpha}$ denotes the root subgroup associated to $\alpha$); if not, then we are done, since then any prime $q_1$ with $\fr_{q_1}=\tilde{\sigma}_1$ (in $\Gal(P_2K_{\phi}K_{\psi}/\tF)$) satisfies the criteria of Hypothesis \ref{hyp1}. But if it is so conjugate, then we must see whether (replacing $\rho_2^Q$ by such a well-placed conjugate) $\kappa(\tilde{\sigma}_1)= \alpha(\rho_2^Q(\tilde{\sigma}_1)) \pmod {\ell^2}$. Again, if these are not equal, we are done. If they are equal, we use the fact that $\Gal(P_2/P_2 \cap K_{\phi} K_{\psi})$ contains $(\Sym^2 \otimes {\det}^{-1})(\bar{r}) \otimes_k \ell W/\ell^2 W$: we can then (using the linear disjointness of $K(\mu_{\ell^2})$ from $P_2 K_{\phi} K_{\psi}$) modify $\tilde{\sigma}_1$ to an element $\sigma_2$ such that $\rho^Q_2(\sigma_2)= \rho^Q_2(\tilde{\sigma}_1) \cdot 2\rho^\vee(1+\ell)$, $\kappa(\sigma_2)= \kappa(\tilde{\sigma}_1)$, $\phi(\sigma_2)=\phi(\tilde{\sigma}_2)$, and $\psi(\sigma_2)=\psi(\tilde{\sigma}_1)$. The existence of primes satisfying Hypothesis \ref{hyp1} follows from the \v{C}ebotarev density theorem.

Having produced $q_1$, to establish the existence of primes $q_2$ satisfying Hypothesis \ref{hyp2} we need to combine the argument of \cite[Theorem 7.4, 10.4]{stp:exceptional} with an analysis of the linear disjointness of the fields $K_{\psi} K_{\phi} K_{\tilde{\phi}}$. The reader may want to review \cite[Lemma 5.3, Theorem 7.4, Lemma 7.6]{stp:exceptional} before proceeding. We first note that $K_{\psi}$ is disjoint over $K$ from $K_{\phi}K_{\tilde{\phi}}$, just as at the end of the proof of \cite[Theorem 7.4]{stp:exceptional}; but in contrast to Proposition \ref{twolindisj}, linear independence of $\phi$ and $\tilde{\phi}$ is not enough to ensure that the fields $K_{\phi}$ and $K_{\tilde{\phi}}$ are linearly disjoint over $K$. Instead, we write $\phi= \sum \phi_i$ and $\tilde{\phi}= \sum \tilde{\phi}_i$, where $\phi_i$ and $\tilde{\phi}_i$ are the respective components of $\phi$ and $\tilde{\phi}$ under the decomposition
\[
H^1(\gal{\Sigma \cup Q_0 \cup q_1}, \br(\fg^\vee)(1))= \bigoplus_i H^1(\gal{\Sigma \cup Q_0 \cup q_1}, \Sym^{2m_i} \otimes {\det}^{-m_i} (\bar{r})(1)).
\]
We then consider the fixed fields $K_{\phi_i}$ and $K_{\tilde{\phi}_j}$ for varying $i$ and $j$. For all $i$, the argument of Proposition \ref{twolindisj} shows that $K_{\phi_i}$ and $K_{\tilde{\phi}_i}$ are either equal or are linearly disjoint over $K$, and that when they are equal we have a linear dependence $\phi_i= \lambda_i \tilde{\phi}_i$ for some $\lambda_i \in k^\times$. There are two cases to consider: suppose first that for all $i$, we have equality, and fix some $i_0$ such that $K_{\phi_{i_0}}= K_{\tilde{\phi}_{i_0}}$ properly contains $K$ (such an $i_0$ exists because the cohomology classes are non-zero). For all $j \neq i_0$, the compositum of the fields $K_{\phi_j}= K_{\tilde{\phi}_j}$ is linearly disjoint from $K_{\phi_{i_0}}$ over $K$ (the different constituents $\Sym^{2m_j} \otimes {\det}^{-m_j} (\bar{r})$ have no common sub-quotient), so we can choose an element $\sigma \in \Gal(K_{\phi}K_{\tilde{\phi}}/K)$ that is trivial in each $\Gal(K_{\phi_j}/K)$, $j \neq i_0$, and takes any desired value in $\Gal(K_{\phi_{i_0}}/K) \xrightarrow{\sim} \Sym^{2m_i} \otimes {\det}^{-m_i} (\bar{r}) (1)$. We can then argue as in \cite[Theorem 7.4, Theorem 10.4]{stp:exceptional}. First, choose $x \in \gal{\tF, \Sigma}$ such that $\br(x)$ is regular semi-simple with $\alpha(\br(x))= \kbar(x)$ for all simple roots $\alpha$, and fix one simple root $\alpha$, requiring in type $\mr{E}_6$ (i.e., $-1 \not \in W_G$) that $\alpha$ is not fixed by the non-trivial pinned outer automorphism. We will modify $x$ by (an extension to $\gal{K}$ of) an element $\sigma \in \Gal(K_{\phi}K_{\tilde{\phi}}/K)$ so that $\phi(\sigma x)$ and $\tilde{\phi}(\sigma x)$ both have non-zero $\mf{g}^\vee_{-\alpha}$ component. To arrange this, we simply take $\sigma$, as above, to be trivial in each $\Gal(K_{\phi_j}/K)$, $j \neq i_0$; and then using \cite[Lemma 7.6]{stp:exceptional} we can choose the restriction of $\sigma$ to $\Gal(K_{\phi_{i_0}}/K)$ to force the desired non-vanishing. Note also that (by \cite[Lemma 7.6]{stp:exceptional}) for such an $\alpha$, and for any non-zero $\psi$, $k[\psi(\gal{K})]$ has non-zero component in the one-dimensional torus generated by the coroot $\alpha^\vee$. This ensures we can make a further modification to $\sigma x$, only affecting its $\Gal(K_\psi/K)$ component, to produce an element $y \in \Gal(K_{\phi}K_{\tilde{\phi}}K_{\psi}/\tF)$ such that the desired primes $q_2$ are those whose frobenius conjugacy classes contain $y$.

The second case is quite similar: suppose now that for some $i$, $K_{\phi_i} \cap K_{\tilde{\phi}_i}= K$. We can then choose $\sigma \in \Gal(K_{\phi} K_{\psi}/K)$ such that $\phi_j(\sigma)= \tilde{\phi}_j(\sigma)=0$ for all $j \neq i$ (by an observation in the previous paragraph, either these values are scalar multiples of each other, or $K_{\phi_j}$ and $K_{\tilde{\phi}_j}$ are linearly disjoint, in which case the values can be arranged independently), and $\phi_i(\sigma)$ and $\tilde{\phi}_i(\sigma)$ assume any desired values. The argument then runs as before.
\endproof
\begin{thm}\label{sl2raising}
Let $G$ be simply-connected of exceptional type, and let $\ell$ be a rational prime greater than $4 h_{G^\vee}-1$, and in the case $G= \mr{E}_8$ not equal to 229, 269, or 367. Let $F$ be a totally real field for which $[F(\zeta_{\ell}):F]= \ell -1$, and let $\bar{r} \colon \gal{F} \to \mr{GL}_2(k)$ be a continuous representation unramified outside a finite set $\Sigma$ of finite places, assumed to contain all places above $\ell$. If $-1$ is not in the Weyl group of $G$, choose a quadratic totally imaginary extension $\tF/F$ in which all elements of $\Sigma$ split, and such that $\tF$ is linearly disjoint from $F(\bar{r}, \zeta_{\ell})$ over $F$; then form the L-group $\LG$ as in \S \ref{ramreview}. Else take $\LG= G^\vee$. Form the composite
\[
\xymatrix{
\gal{\Sigma} \ar[r]^-{\bar{r}} \ar@/^2 pc/[rr]^{\br} & \mr{PGL}_2(k) \times \Gal(\tF/F) \ar[r]^-{\varphi} & {}^L G(k).
} 
\]
Assume that $\bar{r}$ satisfies the following:
\begin{enumerate}
\item For some subfield $k' \subset k$, $\mr{SL}_2(k') \subset \bar{r}(\gal{F}) \subset k^\times \cdot \mr{GL}_2(k')$;
\item $\bar{r}$ is odd;
\item For each $v \vert \ell$, $\bar{r}|_{\gal{F_v}}$ is ordinary, and $\varphi \circ \bar{r}|_{\gal{F_v}}$ satisfies the conditions of Proposition \ref{oldord} below;
\item For each $v \in \Sigma$ not above $\ell$, $\br|_{\gal{\tF_{\tilde{v}}}}$ satisfies a liftable local deformation condition $\mc{P}_v$ whose tangent space has dimension $h^0(\gal{F_v}, \br(\fg^\vee))$. (For $F= \Q$, or more generally $v$ split over $\Q$, we will see in \S \ref{pnotl} that this condition is always satisfied.)
\end{enumerate}
Then there exist a finite set of primes $Q$, disjoint from $\Sigma$ and split in $\tF/F$, and a lift
\[
\xymatrix{ & \LG(\mc{O}) \ar[d] \\
\gal{\Sigma \cup Q} \ar[r]_-{\bar{\rho}} \ar@{-->}[ur]^{\rho^Q} & \LG(k)
}
\]
that has type $\mc{P}_v$ for all $v \in \Sigma$, that has Ramakrishna-type at all $q \in Q$, and that for some $q \in Q$ is ramified.
\end{thm}
\begin{rmk}
The analogue of Proposition \ref{unobstructed} clearly goes through here as well, showing that in the unobstructed case we can still force some Ramakrishna-type ramification. This will be used in Theorem \ref{excapp}.
\end{rmk}
\section{Some local deformation conditions}
\subsection{$p=\ell$}
Let $F/\Ql$ be a finite extension, and suppose $\br \colon \gal{F} \to G^\vee(k)$ factors through a Borel subgroup $B^\vee$ of $G^\vee$. Let $T^\vee= B^\vee/N^\vee$ be the quotient
of $B^\vee$ by its unipotent radical, and fix a lift $\chi \colon I_F \to T^\vee(\mc{O})$ of the composite
\[
 \br|_{I_F} \colon I_F \to B^\vee(k) \to T^\vee(k).
\]
Recall from \cite[Proposition 4.4]{stp:exceptional} that under suitable `regularity' hypotheses the functor of (nearly) ordinary lifts $\Lift_{\br}^{\chi}$ 
(see \cite[Definition 4.1]{stp:exceptional}) satisfies the following:
\begin{prop}[See Proposition 4.4 of \cite{stp:exceptional}]\label{oldord}
 Assume $F$ does not contain $\zeta_{\ell}$, and that $\br$ satisfies the following two conditions:
\begin{itemize}
 \item $H^0(\gal{F}, \br(\fg^\vee/\mf{b}^\vee))=0$. 
\item $H^0(\gal{F}, \br(\fg^\vee/\mf{b}^\vee)(1))=0$.
\end{itemize}
Then $\Lift^{\chi}_{\br}$ is a liftable local deformation condition with tangent space $L_{\br}^{\chi}$ of dimension 
\[
[F:\Ql]\dim(G^\vee/B^\vee)+ h^0(\gal{F}, \br(\fg^\vee).
\]
\end{prop}
In the application in \cite{stp:exceptional}, these two vanishing hypotheses were arranged by using $\br$ such that for all simple roots $\alpha$, $\alpha \circ \br|_{I_F}$
was equal to $\kbar^{r_{\alpha}}$ for some integers $r_{\alpha} \geq 2$. Namely, the $\br$ in question were constructed by taking a classical modular form $f$ of weight at least 3,
with associated $r_{f, \ell} \colon \gal{\Q} \to \mr{PGL}_2(\overline{\Z}_{\ell})$, and setting $\br= \varphi \circ \bar{r}_{f, \ell}$. We will now check that the two regularity hypotheses in fact
still hold, for a density 1 set of $\ell$, for the Galois representations associated to elliptic curves; this is a consequence of work of Weston,\footnote{See \cite[Proposition 4.4, Theorem 5.5]{weston:unobstructed}. Note, however, that his final claim is only valid--or at least proven--for elliptic curves, not general weight two modular forms.}
to which Shekhar Khare drew my attention. 
\begin{prop}[See Proposition 4.4 and Theorem 5.5 of \cite{weston:unobstructed}]\label{ordweight2} 
Let $E/\Q$ be an elliptic curve. Then there is a density one set of rational primes $\ell$, with associated (projective) $\ell$-adic representation $r_{E, \ell} \colon \gal{\Q} \to \mr{PGL}_2(\Z_{\ell})$, such that the composite $\br= \varphi \circ \bar{r}_{f, \ell}$ is ordinary and satisfies
\begin{itemize}
 \item $H^0(\gal{\Ql}, \br(\fg^\vee/\mf{b}^\vee))=0$; and
\item $H^0(\gal{\Ql}, \br(\fg^\vee/\mf{b}^\vee)(1))=0$.
\end{itemize}
Indeed, for these conditions to hold, it suffices that $\ell \geq h_{G^\vee}$ and $a_{\ell}^2 \not \equiv 1 \pmod \ell$.
\end{prop}
\proof
By considering the $\gal{\Ql}$-action on $\br(\fg^\vee/\fb^\vee)$ with a basis of root spaces ordered by root height, we see that he first vanishing condition holds even on $I_{\Ql}$ as long as $\ell \geq h_{G^\vee}$, so it suffices to attend to the second vanishing condition; and for the same reason of the inertial action, we need only show that there are no invariants on (the image modulo $\fb^\vee$ of) the negative simple root spaces. Recall (\cite[Theorem 2]{wiles:ordinary}) that if $\ell$ is ordinary for $E$, then $r_{E, \ell}|_{\gal{\Ql}}$ has the form 
\[
r_{E, \ell}|_{\gal{\Ql}} \cong \begin{pmatrix}
\kappa \cdot \lambda_{\alpha^{-1}} & * \\
0 & \lambda_{\alpha} \\
\end{pmatrix},
\]  
where $\lambda_\beta$ denote the unramified character by which (arithmetic) Frobenius acts by $\beta$, and where $\alpha \in \Z_{\ell}$ is the unit root of the polynomial $X^2-a_{\ell}X+\ell=0$. A quick calculation shows that $h^0(\gal{\Q_{\ell}}, \br(\fg^\vee/\fb^\vee)(1))=0$ as long as $\alpha^2 \not \equiv 1 \pmod \ell$, or equivalently $a_{\ell}^2 \not \equiv 1 \pmod \ell$. Since $a_{\ell}$ is an integer of absolute value at most $2 \sqrt{\ell}$, $a_{\ell} \equiv \pm 1 \pmod \ell$ forces (for $\ell \geq 7$) $a_{\ell}= \pm 1$. To ensure $\ell$ was ordinary for $E$ in the first place, we needed (for $\ell >>0$) $a_{\ell} \neq 0$, so in sum we can apply \cite[Theorem 20]{serre:cebotarev} to find a density 1 set of places $\ell$ at which $a_{\ell}$ avoids the values $\{-1, 0, 1\}$, completing the proof of the Proposition.
\endproof

\subsection{$p \neq \ell$}\label{pnotl}
In the construction of Galois representations with exceptional monodromy groups in this paper, we will deform Galois representations of the form $\br= \varphi \circ \bar{r}_{E, \ell} \colon \gal{\Q} \to \LG(k)$, where $E$ is an elliptic curve over $\Q$. We will have specified the local behavior of $E$ at two finite places ($\ell$ and one auxiliary semi-stable prime), but otherwise we will have no (effective) control over the bad reduction of $E$ (see \S \ref{exceptional}). Thus, we will need to be able to define locally liftable deformation conditions with big enough tangent space for all primes, and in this section we will construct `minimal' deformation conditions for residual representations of the form $\varphi \circ \bar{r} \colon \gal{\Q_p} \to G^\vee(k)$, starting from some $\bar{r} \colon \gal{\Q_p} \to \mr{GL}_2(k)$. This is rather \textit{ad hoc}, but it suffices for our purposes; Jeremy Booher has understood more systematically how to construct minimal deformation conditions in a general setting (\cite{booher:thesis}). Throughout this section we assume $p \neq \ell$, and for simplicity we assume $\ell \geq 5$ (in the application of \S \ref{exceptional}, we will make more restrictive hypotheses on $\ell$). We recall also that we always assume $G^\vee$ is an adjoint group, so $\br$ factors through the projectivization $\mathbb{P}(\bar{r})$ of $\bar{r}$ (this hypothesis can certainly be removed, if desired).

First, recall Diamond's classification of two-dimensional representations $\bar{r} \colon \gal{\Q_p} \to \mr{GL}_2(k)$ into four cases: principal, special, vexing, and harmless (\cite[\S 2]{diamond:extension}). Because $\ell \geq 5$, the image $\bbP(\bar{r})(I_{\Q_p})$ has order prime to $\ell$ unless we are in the `special' case (see \cite[Proposition 2.1, 2.3, 2.4]{diamond:extension}), so we immediately deduce the following lemma:
\begin{lemma}\label{min1}
Suppose $p \neq \ell$ and $\ell \geq 5$. If $\bar{r}$ is not special, then the minimal deformation condition of \cite[\S 4.4]{stp:exceptional} is a liftable deformation condition for $\br$ with tangent space of dimension $h^0(\gal{\Q_p}, \br(\fg^\vee))$.
\end{lemma}
It thus only remains to treat the case when $\bar{r}$ is special. Then by \cite[Proposition 2.2]{diamond:extension}, $\bbP(\bar{r})(I_{\Q_p})$ is cyclic of order $\ell$, and $\bbP(\bar{r})$ factors through the tame quotient $T_{\ell}= \Gal(\Q_p^{\mr{ur}}(p^{1/\ell^{\infty}})/\Q_p)$, and we will restrict to representations of this group, which is isomorphic to the semi-direct product $\Z_{\ell} \rtimes \hat{\Z}= \langle \tau \rangle \rtimes \langle \Fr_p \rangle$, with the familiar action $\Fr_p \tau \Fr_p^{-1}= \tau^p$. Denoting by $B^\vee \supset N^\vee$ the Borel (and unipotent radical) used to define our principal homomorphism $\varphi$, we can and do reduce to the case where $\bar{r}(\tau)= \begin{pmatrix} 1 & 1 \\ 0 & 1 \end{pmatrix}$, so $\br(\tau)$ is a regular unipotent element in $N^\vee(k)$.
\begin{defn}
Let $\Lift^{\mr{min}}_{\br}$ be the sub-functor of $\Lift_{\br}$ consisting of all lifts $\rho \colon \gal{\Q_p} \to G^\vee(R)$ such that 
\begin{itemize}
\item $\rho$ factors through $T_{\ell}$; and
\item $\rho(\tau)$ is $\widehat{G^\vee}(R)$-conjugate to an element of $N^\vee(R)$.
\end{itemize}
\end{defn}
$\Lift^{\min}_{\br}$ is obviously a sub-functor of $\Lift_{\br}$.
\begin{lemma}\label{min2}
Assume $\ell$ is very good for $G^\vee$ (eg, $\ell > h_{G^\vee}$ suffices). $\Lift^{\mr{min}}_{\br}$ is a liftable local deformation condition with (after passing to deformation classes) tangent space equal to the unramified cohomology $H^1(\gal{\Q_p}/I_{\Q_p}, \br(\fg^\vee)^{I_{\Q_p}})$, hence of dimension $h^0(\gal{\Q_p}, \br(\fg^\vee))$.
\end{lemma}
\proof
The key observation, used in each step of the proof, is that if $\gamma \in N^\vee(k)$ is regular unipotent, then 
\begin{equation}\label{key}\tag{REG}
1- \Ad(\gamma) \colon \fb^\vee \to \mf{n}^\vee
\end{equation}
is surjective, and that the kernel of $1- \Ad(\gamma)$ on all of $\fg^\vee$ is in fact contained in $\mf{n}^\vee$. Surjectivity follows by a dimension-count: since $\ell$ is very good for $G^\vee$, the centralizer of $\gamma$ in $\fg^\vee$ is the Lie algebra of the centralizer in $G^\vee$, of dimension $\rk \fg^\vee$.

First, the Mayer-Vietoris property is checked by an argument similar to that of \cite[Lemma 4.10]{stp:exceptional}, the key point being that if a lift $\rho$ has the property that $\rho(\tau) \in N^\vee(R)$ and that $g \rho(\tau) g^{-1} \in N^\vee(R)$ for some $g \in G^\vee(R)$, then $g \in B^\vee(R)$ (an induction argument using Equation (\ref{key})). We omit the details.

An easy calculation in $\mr{GL}_2(k)$ with the relation $\Fr_p \tau \Fr_p^{-1}= \tau^p$ shows that $\br(\Fr_p)$ is an element of $B^\vee(k)$. Let $R \to R/I$ be a small surjection. Up to $\widehat{G^\vee}(R/I)$-conjugation, an element $\rho$ of $\Lift^{\mr{min}}_{\br}(R/I)$ corresponds to a pair $t= \rho(\tau) \in N^\vee(R/I)$ and $g= \rho(\Fr_p)$ satisfying $gtg^{-1}=t^p$. By the previous paragraph, $g$ belongs to $B^\vee(R/I)$. Choose any lifts $\tilde{g}$ and $\tilde{t}$ of $g$ and $t$ to $B^\vee(R)$ and $N^\vee(R)$; then 
\[
\exp(X)= \tilde{g}\tilde{t}\tilde{g}^{-1} \tilde{t}^{-p}
\]
for some $X \in \mf{n}^\vee \otimes_k I$. We claim that, retaining the choice of $\tilde{t}$, we can modify $\tilde{g}$ to an element $\exp(Y) \tilde{g}$, for some $Y \in \fg^\vee \otimes_k I$, such that $(\Fr_p, \tau) \mapsto (\exp(Y) \tilde{g}, \tilde{t})$ corresponds to an element of $\Lift^{\mr{min}}_{\br}(R)$ lifting $\rho$. It is equivalent to check that we can choose $Y$ such that $X+Y- \Ad(\br(\tau)^{-p})Y=0$. By our key observation (Equation (\ref{key})), the linear map
\[
1- \Ad(\br(\tau)^{-p}) \colon \fb^\vee \to \mf{n}^\vee
\]
is surjective since $\br(\tau)^{-p}$ is a regular element. 

Finally, we compute the tangent space. Since the invariants $\br(\fg^\vee)^{I_{\Q_p}}$ equal the centralizer in $\fg^\vee$ (hence in $\mf{n}^\vee$) of the regular unipotent element $\br(\tau)$, it is clear that the tangent space $L^{\mr{min}}_{\br}$ contains the unramified cohomology $H^1(\gal{\Q_p}/I_{\Q_p}, \br(\fg^\vee)^{I_{\Q_p}})$. But again the surjectivity of $1- \Ad \br(\tau) \colon \fb^\vee \to \mf{n}^\vee$ allows us to trivialize (i.e., realize as a coboundary) the restriction to $I_{\Q_p}$ of any cocycle $\phi$ corresponding to an element of $\Lift^{\mr{min}}_{\br}(k[\epsilon])$. Thus, $H^1(\gal{\Q_p}/I_{\Q_p}, \br(\fg^\vee)^{I_{\Q_p}})$ is in fact the entire tangent space of $\Def^{\mr{min}}_{\br}$.
\endproof
\section{Application to exceptional monodromy}\label{exceptional}
In this section we will make two improvements to the construction of geometric Galois representations with exceptional monodromy groups given in \cite{stp:exceptional}. One could hope to improve (at least) two aspects of the results of \cite{stp:exceptional}:
\begin{itemize}
\item in \cite{stp:exceptional}, geometric representations $\rho_{\ell} \colon \gal{\Q} \to \LG(\Qlb)$ with Zariski-dense image are constructed only for a density one set of rational primes $\ell$;
\item the $\rho_{\ell}$ above are constructed to have `generic' Hodge-Tate weights, something perfectly natural from the point of view of $\ell$-adic deformation, but rather undesirable from a `motivic' perspective: in particular, because of Griffiths' transversality, such $\rho_{\ell}$ cannot appear as specializations of an arithmetic local system over a curve (in contrast to the $\mr{G}_2$, $\mr{E}_7$, and $\mr{E}_8$ examples of \cite{yun:exceptional}.
\end{itemize}
The basic idea for improving the `density-one' restriction is that rather than starting from a single well-chosen modular form, and letting $\ell$ range over its (density one) set of ordinary primes, we instead for each $\ell$ choose a different elliptic curve to provide the `seed' two-dimensional $\gal{\Q}$-representation. It seems to be more difficult, given $\ell$, to produce a weight-three modular form that serves our purpose, so we are stuck working in weight two and therefore need Proposition \ref{ordweight2}. The elliptic curve we construct will have some specific local behavior at three different primes, but we will not be able to control its local behavior at other primes (this seems to be related to quite deep questions in analytic number theory); thus we require the local results of \S \ref{pnotl}. 

To achieve big-monodromy lifts whose Hodge numbers are ``consecutive'' (see Theorem \ref{excapp}), we use the level-raising results of \S \ref{sl2section}: those techniques allow us to construct lifts whose image contains a principal $\mr{SL}_2$ \textit{and} an element of some non-regular (and non-trivial) unipotent conjugacy class, which suffices to ensure full monodromy except in the case $G= \mr{E}_6$ (where, unfortunately, I don't see how to improve the argument). Here then is the improvement made possible by the results of the preceding sections: 
\begin{thm}\label{excapp}
Let $G$ be a simply-connected group of exceptional type, with the exception of $\mr{E}_6$. Let $\ell$ be any rational prime greater than $4h_{G^\vee}-1$, and also exclude $\ell= 229, 269, 367$ for $G$ of type $\mr{E}_8$. Then there are infinitely many (non-conjugate) geometric representations $\rho_{\ell} \colon \gal{\Q} \to \LG(\Ql)$ with Zariski-dense image whose Hodge-Tate co-character is (up to conjugacy) equal to $\rho^\vee \in X_{\bullet}(T^\vee)$; more generally, the same holds for any other co-character $\mu \in X_{\bullet}(T^\vee)$ satisfying $\alpha \circ \mu >0$ and $\alpha \circ \mu \equiv 1 \pmod {\ell -1}$ for all simple roots $\alpha$ (here we identify $\Hom(\mbf{G}_m, \mbf{G}_m) = \Z$, mapping the identity to 1).

If $G$ is of type $\mr{E}_6$, the same assertion holds but without the assertion about the Hodge-Tate co-character.
\end{thm}
\proof
Fix a prime $\ell > 4h_{G^\vee} -1$, and not equal to 229, 269, or 367 if $G= \mr{E}_8$. We will show there exists an elliptic curve over $\Q$ with the following local behavior:
\begin{itemize}
\item $E$ is (good) ordinary at $\ell$ with Hecke polynomial $X^2-a_{\ell} X + \ell$ for an integer $a_{\ell}$ not congruent to 0 or $\pm 1$ modulo $\ell$.
\item For some auxiliary prime $p_0$ such that the order of $p_0$ modulo $\ell$ is at least $h_{G^\vee}$, $E \otimes_{\Q} \Q_{p_0}$ has split multiplicative reduction, and the $\gal{\Q_{p_0}}$-action on $E[\ell](\overline{\Q}_{p_0})$ is non-split. 
\item For some auxiliary prime $p_1$, $E \otimes_{\Q} \Q_{p_1}$ has additive reduction, and the $\gal{\Q_{p_1}}$-action on $E[\ell](\overline{\Q}_{p_1})$ is irreducible.
\item At all primes $p \not \in \{p_0, p_1, \ell\}$, we impose no restriction on the ramification of $E \otimes_{\Q} \Q_p$.
\end{itemize}
Here is how to `produce' such an $E/\Q$. Choose an integer $a$ such that $(a, \ell)=1$, $|a| < 2 \sqrt{\ell}$, and $a \neq \pm 1$. By Honda-Tate theory, there is an elliptic curve over $\mathbb{F}_{\ell}$ corresponding to the Weil $\ell$-number roots of $X^2-aX+\ell=0$. Fix any Weierstrass equation over $\mathbb{F}_{\ell}$ for this curve. Next, take any prime $p_0$ whose order mod $\ell$ is at least $h_{G^\vee}$, and consider the Weierstrass equation 
\[
y^2+xy= x^3-5p_0 x-p_0 \pmod {p_0^2}.
\]
This is the reduction modulo $p_0^2$ of a Tate curve over $\Z_{p_0}$ with Tate parameter $p_0$, and in particular having non-split mod $\ell$ representation (see \cite[V Theorem 3.1, Proposition 6.1]{silverman:aec2}; any lift of this Weierstrass equation is also a Tate curve with Tate parameter having $p_0$-adic valuation equal to 1 (\cite[V Theorem 5.3]{silverman:aec2}). Finally, let $p_1 \geq 5$ be a prime congruent to 2 mod 3, and consider the additive reduction Weierstrass equation
\[
y^2=x^3-3p_1x-2p_1 \pmod {p_1^2}.
\]
The same equation regarded over $\Q_{p_1}$ has $j$-invariant determined by $p_1= \frac{j}{j-1728}$ and discriminant $\Delta= \frac{12^6 j^2}{(j-1728)^3}$ (see \cite[Equation 0.2]{diamond-kramer}), so $\ord_{p_1}(j)=1$, $\ord_{p_1}(j-1728)=0$, and $\ord_{p_1}(\Delta)=2$. Then by \cite[Table 1, Proposition 0.3]{diamond-kramer} (note that our hypothesis implies $\mu_3$ is not contained in $\Q_{p_1}$), for any lift to $\Z_{p_1}$ of this $\pmod {p_1^2}$ equation, the associated elliptic curve over $\Q_{p_1}$ will have irreducible (`type V') mod $\ell$-representation.

We now have three Weierstrass equations, modulo $\ell$, $p_0^2$, and $p_1^2$. We solve for a common lift to a Weierstrass equation with integral coefficients, and we let $E$ be the associated elliptic curve over $\Q$. First we claim that $\bar{r}_{E, \ell} \colon \gal{\Q} \to \mr{GL}_2(\mathbb{F}_{\ell})$ is surjective. Since $\bar{r}_{E, \ell}|_{\gal{\Q_{p_0}}}$ is non-split, the image contains an element of order $\ell$; by \cite[2.4 Proposition 15]{serre:propgal}, the image of $\bar{r}_{E, \ell}$ is then either contained in a Borel, or contains $\mr{SL}_2(\mathbb{F}_{\ell})$. But irreducibility of $\bar{r}_{E, \ell}|_{\gal{\Q_{p_1}}}$ rules out the Borel case, and since $\det(r_{E, \ell})$ is surjective, we conclude that the image of $\bar{r}_{E, \ell}$ is all of $\mr{GL}_2(\mathbb{F}_{\ell})$. Let $\br \colon \gal{\Q} \to \LG(\mathbb{F}_{\ell})$ be the composite $\varphi \circ \bar{r}_{E, \ell}$. The restriction $\br|_{\gal{\Ql}}$ is ordinary and satisfies the conclusion of Proposition \ref{ordweight2}; here we use $a_{\ell}^2 \not \equiv 1 \pmod \ell$. It follows (Proposition \ref{oldord}) that the ordinary lift functor $\Lift_{\br}^{\chi}$ associated to the character $\chi \colon I_{\Ql} \to T^\vee(\mathbb{Z}_{\ell})$ determined by $\alpha \circ \chi= \kappa$ for all simple roots $\alpha$ is formally smooth with tangent space of the correct dimension. At the prime $p_0$, we take a Steinberg deformation condition as in \cite[4.3]{stp:exceptional}. At all primes $p \not \in \{\ell, p_0\}$, we take (liftable) minimal deformation conditions as in Lemmas \ref{min1} or \ref{min2}. We can therefore apply \cite[Theorem 7.4, Theorem 10.4]{stp:exceptional} to produce geometric lifts of $\br$ with the local behavior just indicated, and it follows as in \cite[Theorem 8.4]{stp:exceptional} that the monodromy group of such a lift 
\[
\rho^Q \colon \gal{\Q} \to \LG(\Z_{\ell})
\]
($Q$ here indicates an auxiliary set of Ramakrishna primes, disjoint from $\ell$ and the bad primes of $E$) contains the principal $\mr{SL}_2$.

Now we combine Proposition \ref{forceramsl2} and (the proofs of) Propositions \ref{formalforcing} and \ref{unobstructed} to produce a new lift $\rho^{\widetilde{Q}} \colon \gal{\Q} \to \LG(\Z_{\ell})$ for some auxiliary set $\widetilde{Q}$, such that:
\begin{itemize}
\item at all primes not in $\widetilde{Q}$, $\rho^{\widetilde{Q}}$ has the same local behavior just indicated;
\item at at least one prime in $\widetilde{Q}$, $\rho^{\widetilde{Q}}$ is ramified (of Ramakrishna type).
\end{itemize}
Note that in the principal $\mr{SL}_2$ case (Proposition \ref{forceramsl2}), we have only carried out the construction of ramification-forcing auxiliary primes for the mod $\ell^2$ representation; this has the effect that we do not guarantee $\rho^{\widetilde{Q}}$ is ramified at all places in $Q$ (the argument does not retain the ramification at the set denoted $Q_{\mr{ram}}$ in the proof of Proposition \ref{formalforcing}). But no matter: having any Ramakrishna-type ramification implies that the monodromy group of $\rho^{\widetilde{Q}}$ cannot be contained in the principal $\mr{SL}_2$. To see this, let $q$ be a prime in $\widetilde{Q}$ at which $\rho^{\widetilde{Q}}$ is in fact ramified, and let $t_{q, \ell}$ be a generator of the $\ell$-part of the tame inertia group $I_{\Q_q}^t$. Certainly $\rho^{\widetilde{Q}}|_{I_{\Q_q}}$ factors through $I_{\Q_q}^t$, and we claim that $\rho^{\widetilde{Q}}(t_{q, \ell})$ is unipotent. Up to conjugation, we may assume $\rho^{\widetilde{Q}}(\gal{\Q_q})$ lies in the semi-direct product $U_{\alpha^\vee} \rtimes T^\vee$; we write it as a pair $(u, s)$ with $u \in U_{\alpha^\vee}(\Z_{\ell})$, $s \in T^\vee(\Z_{\ell})$. The conjugation relation $fr_{q}^{-1} t_{q, \ell} fr_q= t_{q, \ell}^q$ implies $s=s^q$, but combining the fact that $s \pmod \ell=1$ and $q \not \equiv 1 \pmod \ell$ forces $s=1$ (one can argue this by induction, looking at the images $s \pmod {\ell^n}$ for all $n \geq 1$). 

The (non-trivial) unipotent conjugacy class of $\rho^{\widetilde{Q}}(t_{q, \ell})$ is not the principal orbit, so the monodromy group of $\rho^{\widetilde{Q}}$ cannot be a principal $\mr{SL}_2$. Therefore by Dynkin's theorem the monodromy group of $\rho^{\widetilde{Q}}$ must either be all of $\LG$ or, in the case $G= \mr{E}_6$, a copy of $\mr{F}_4$.

In the case $G= \mr{E}_6$, we repeat the same argument, but instead of taking ordinary deformations associated to the character $\chi \colon I_{\Q_{\ell}} \to T^\vee(\Z_{\ell})$ satisfying $\alpha \circ \chi= \kappa$ for all $\alpha \in \Delta$, we choose $\chi$, as in \cite[Theorem 8.4, Theorem 10.6]{stp:exceptional}, satisfying $\alpha \circ \chi= \kappa^{r_{\alpha}}$ for distinct positive integers $r_{\alpha} \equiv 1 \pmod {\ell -1}$.\footnote{We note there is a typo in \cite[Theorem 8.4]{stp:exceptional}: the integers $r_{\alpha}$ there must satisfy $r_{\alpha} \equiv 2 \pmod {\ell -1}$, since, in the notation of that argument, $\alpha \circ \bar{\rho}_T= \bar{\kappa}^2$.}
\endproof
\begin{rmk}
The prime $p_1$ is used to make sure the image of $\bar{r}_{E, \ell}$ is not contained in a Borel. A heavy-handed alternative to using $p_1$ would be to invoke Mazur's theorem on rational isogenies, which, given the local behavior at $p_0$, implies (\cite[Theorem 3]{mazur:rational}) surjectivity of $\bar{r}_{E, \ell}$ as long as we avoid the set $\ell \leq 19$ and $\ell \in \{37, 43, 67, 163\}$; note that most of these possibilities are already excluded by the assumption $\ell > 4h_{G^\vee} -1$.
\end{rmk}
\bibliographystyle{amsalpha}
\bibliography{biblio.bib}

\end{document}